\input amstex
\documentstyle{amsppt}
\input epsf
\magnification 1200
\NoBlackBoxes
\vcorrection{-9mm}

\topmatter
\title   Two-dimensional diffusion orthogonal polynomials
         ordered by a weighted degree
\endtitle
\author   S.~Yu.~Orevkov
\endauthor

\def\bR{\Bbb R}
\def\cL{\Cal L}
\def\LL{\bold{L}}

\abstract
We study the following problem: describe the triplets $(\Omega,g,\mu)$,
$\mu=\rho\,dx$, where $g= (g^{ij}(x))$ is the (co)metric associated
with the symmetric second order differential operator
$\LL (f) = \frac{1}{\rho}\sum_{ij} \partial_i (g^{ij} \rho \partial_j f)$ 
defined on a domain $\Omega$ of $\bR^d$ and such that there exists an
orthonormal basis  of $\cL^2(\mu)$ made of polynomials which are eigenvectors
of $\LL$, where the polynomials are ranked according to some weighted degree.

In a joint paper with D.~Bakry and
M.~Zani this problem was solved in dimension 2 for the usual degree.
In the present paper we solve it still in dimension 2, but
for a weighted degree with arbitrary positive weights.
\endabstract
\address
Steklov Mathematical Institute, Gubkina 8, Moscow, Russia
\endaddress

\address
IMT, l'universit\'e Paul Sabatier, 118 route de Narbonne, Toulouse, France
\endaddress

%\address
%AGHA Laboratory, Moscow Institute of Physics and Technology, Russia
%\endaddress

\email
orevkov\@math.ups-tlse.fr
\endemail

\endtopmatter
\rightheadtext{Two-dimensional diffusion orthogonal polynomials}

%========================= References ===============
\def\refBB      {1}
\def\refBZ      {2}
\def\refBOZ     {3}
\def\refKu      {4}
\def\refSouArx  {5}
\def\refSouAFST {6}

%======================= References to [\refBOZ] ============
\def\BOZeqChangeVar   {2.3}
\def\BOZpropChangeVar {2.5}
\def\BOZpropSDOP      {2.12}

%============================ Sections ====================

\def\remDimOne  {1.1}
\def\defDOP     {1.2}

\def\sectGen {2}
   \def\sectAlgDOP    {\sectGen.1}
   \def\sectChangeVar {\sectGen.2}
   \def\sectALGtoDOP  {\sectGen.3}

\def\sectDimTwoAlg {3}
   \def\sectW        {\sectDimTwoAlg.2}
   \def\sectLocal    {\sectDimTwoAlg.3}

\def\sectBigW      {4}

\def\sectSmallW    {5}
   \def\sectHirz     {\sectSmallW.1}
   \def\sectIrred    {\sectSmallW.2}
   \def\sectDegTwo   {\sectSmallW.3}

\def\sectSDOP      {6}
   \def\sectCompact  {\sectSDOP.1}
   \def\sectNonComp  {\sectSDOP.2}

\def\sectRealiz    {7}

%=========================== Proclaims ====================
\def\propAlgDOP    {2.1}
\def\defAlgDOP     {2.2}
\def\propW         {2.3}
\def\propGammaOne  {2.4}

\def\propChangeVar  {2.5}
\def\exaChangeVar   {2.6}
\def\exaChangeVarII {2.7}

\def\propOPlus     {2.8}
\def\propDOP       {2.9}
\def\defMaxBndry   {2.10}
\def\propMes       {2.11}
\def\remMes        {2.12}
\def\remMesBis     {2.13}
\def\corMes        {2.14}
\def\propMesII     {2.15}

%============================
\def\lemNewton    {3.1}
\def\lemChangeW   {3.2}
\def\lemValuA     {3.3}
\def\lemValu      {3.4}

%============================
\def\lemA         {4.1} %% no branch at (0,infty)
\def\lemSym       {4.2}
\def\propDegGaTwo {4.3}
\def\propDegGaOne {4.4}

%============================
\def\lemPlu       {5.1}
\def\lemParam     {5.2}
\def\propIrred    {5.3}
\def\remIrred     {5.4}
\def\propReducTwo {5.5}
\def\lemReduc     {5.6}
\def\propReducThree {5.7}

%============================
\def\thCompTwo  {6.1}
\def\thCompInf  {6.2}
\def\thNonCTwo  {6.3}
\def\thNonCInf  {6.4}
\def\remCorrig  {6.5}

%============================= Tables =====================

\def\tabQuo   {1}

%============================ Equations ===================
\def\eqLL       {1}
\def\eqSym      {2}
\def\eqLLrev    {3}
\def\eqAlgDOP   {4}
\def\eqChangeVar{5}
\def\eqDh       {6}
\def\eqDkDj     {7}
\def\eqMes      {8}
\def\eqMesIneq  {9}
\def\eqBord     {10}
\def\eqParam    {11}
\def\eqValu     {12}
\def\eqPfLemA   {13}
\def\eqGDtwoA   {14}
\def\eqGDtwoB   {15}
\def\eqGDoneA   {16}
\def\eqLaurent  {17}
\def\eqGDoneB   {18}
\def\eqWeiPro   {19}
\def\eqHirz     {20}

\def\eqGenusP   {21}
\def\eqGenus    {22}
\def\eqRH       {23}
\def\eqDual     {24}
\def\eqI        {25}
\def\eqBlowUp   {26}
\def\eqIII      {27}
\def\eqIG       {28}
\def\eqReduc    {29}
\def\eqReducIII {30}
\def\eqCompInf  {31}
\def\eqL        {32}
\def\eqDhh      {33}
\def\eqM        {34}

%===================== AlgDOP conditions =================

\def\condADi   {A1}
\def\condADii  {A2}
\def\condADiii {A3}

%======================== Model numbers ==================

\def\modelBi   {B1}   % Bounded
\def\modelBii  {B2}
\def\modelBiii {B3}
\def\modelBiv  {B4}
\def\modelBv   {B5}

\def\modelUi   {U1}
\def\modelSum  {U2}

%========================== Figures ======================
\def\figABCi   {1}
\def\figABCii  {2}
\def\figHirz   {3}
\def\figBlowUp {4}
\def\figAA     {5}
\def\figIV     {6}
\def\figIVtan  {7}
\def\figCompTwo{8}
%=========================  Macros ======================
\def\bR{\Bbb{R}}
\def\bC{\Bbb{C}}

\def\CP{\Bbb{CP}}
\def\bK{\Bbb{K}}
\def\bP{\Bbb{P}}
\def\bN{\Bbb{N}}
\def\bZ{\Bbb{Z}}

\def\bS{\Bbb{S}}
\def\bB{\Bbb{B}}
\def\Re{\operatorname{Re}}
\def\ord{\operatorname{ord}}
\def\diag{\operatorname{diag}}
\def\const{\operatorname{const}}

\def\LL{\bold{L}}
\def\cL{\Cal{L}}
\def\cF{\Cal{F}}
\def\cP{\Cal{P}}
\def\cN{\Cal{N}}
\def\ww{\bold w}
\def\xx{\bold x}

\def\diff#1#2{#1'_#2}
\def\g{g}
\loadbold
\def\De{{\boldsymbol\Delta}}

\document

\head 1. Introduction
\endhead

In [\refBOZ] the following problem is studied
(see also [\refBB], [\refBZ], [\refSouArx], [\refSouAFST]): describe all triples
$(\Omega,\LL,\mu)$ where $\Omega$ is a domain in $\bR^d$, $\LL$ is
an elliptic second order operator of the form
$$
      \LL(f) = \sum_{i,j} g^{ij}(x)\partial_{ij}f
             + \sum_ib^i(x)\partial_i f                                 \eqno(\eqLL)
$$
with $g^{ij}$ and $b^i$ continuous
in $\Omega$ and $\mu$ a probability measure on $\Omega$
which is absolutely continuous with respect to the
Lebesgue measure, and such that there exists a polynomial orthogonal basis of
$\cL^2(\Omega,\mu)$ formed by eigenvectors of $\LL$, which is also a basis
(in the algebraic sense) of $\bR[x]$, $x=(x_1,\dots,x_d)$. It is clear that
in this case $\LL$ is symmetric on the space of polynomials, i.e.,
$$
   \int_\Omega f_1\LL f_2\,d\mu = \int_\Omega f_2\LL f_1\,d\mu          \eqno(\eqSym)
$$
for any two polynomial functions $f_1$ and $f_2$.

Suppose that $e_1,e_2,\dots$ is such a basis, and let $V_n$ be the subspace spanned
by $e_1,\dots,e_n$. Then we have an increasing sequence of $\LL$-invariant
subspaces $V_1\subset V_2\subset V_3\subset\dots$ of $\bR[x]$ 
whose union is $\bR[x]$, and $\bR[x]$ is dense in $\cL^2(\Omega,\mu)$.

Vice versa, given an increasing sequence of finite-dimensional
$\LL$-invariant subspaces of $\bR[x]$ whose union is $\bR[x]$, one can always choose
an orthogonal polynomial eigenbasis of $\cL^2(\Omega,\mu)$ provided
that $\LL$ is symmetric on polynomials and the space of polynomials
is dense in $\cL^2(\Omega,\mu)$.

It does not seem that this problem can be solved in such generality,
without imposing a condition that the filtration $V_1\subset V_2\subset\dots$
is somewhat natural. So, in [\refBOZ] the above problem was considered with
an additional condition that for any $n$, the space of polynomials
of degree $\le n$ is invariant under $\LL$ and thus occurs among the $V_i$'s.
It was also assumed in [\refBOZ] that,
when $\Omega$ is not bounded, (\eqSym) holds for any pair of
compactly supported functions (for bounded domains the latter condition easily
follows from the symmetry of $\LL$ on the polynomials combined with the density
of $\bR[x]$ in $\cL(\Omega,\mu)$). Under these assumptions,
a complete list of solutions of the above problem is given in [\refBOZ]
in dimensions $2$. Up to affine transformation of $\bR^2$, there is a
one-parameter family of bounded domains and also 18 rigid domains
(10 of them are bounded) for which there exists a solution.

\medskip\noindent
{\bf Remark \remDimOne.}
In dimension 1, the only solutions under the aforementioned assumptions are
the classical systems of orthogonal polynomials: Hermite, Laguerre, and Jacobi
polynomials. They are obtained by Gram-Schmidt orthogonalization
process for the measure densities, respectively,
$Ce^{-x^2/2}$ on $\bR$, $C_ax^{a-1}e^{-x}$ with $a>0$ on $[0,\infty)$, and
$C_{p,q}(1-x)^{p-1}(1+x)^{q-1}$ with $p,q>0$ on $[-1,1]$ (here $C$, $C_{a}$, and
$C_{p,q}$ are normalizing constants).
The corresponding operators are $\partial^2-x\partial$, $x\partial^2+(a-x)\partial$,
and $J_{p,q}=(1-x^2)\partial^2-\big((p-q)+(p+q)x\big)\partial$.
\medskip

It turns out that the filtration by the usual degree is too restrictive in
dimension $d\ge2$. Several
natural systems of orthogonal polynomial are not covered. However,
they can be obtained by this procedure if one considers a weighted degree
instead (see [\refBB]). As usually, the weighted degree of a polynomial
$P=\sum_k a_k x^k$, $k=(k_1,\dots,k_d)$, with real positive weights
$\ww=(w_1,\dots,w_d)$ is defined as
$$
    \deg_{\ww}(P)=\max_{a_k\ne 0}(w_1k_1+\dots+w_dk_d).
$$
In this paper (in \S\sectSDOP),
for any pair of positive weights $\ww=(w_1,w_2)$, we give a complete list
of two-dimensional solutions satisfying the condition that for any $n$, the set of
polynomials $P$ with $\deg_{\ww}(P)\le n$ is invariant under $\LL$.

Let us give precise definitions. We say that $\Omega\subset\bR^d$
is a {\it natural domain} if it is a connected open set which coincides
with the interior of its closure.

\medskip\noindent
{\bf Definition \defDOP.} (cf.~[\refBB, \refBOZ])
Let $\Omega\subset\bR^d$ be a natural domain, $\LL$ be an elliptic second
order differential operator of the form (\eqLL) with coefficients continuous
in $\Omega$, and $\mu(dx)=\rho(x)\,dx$ be a probability measure on $\Omega$
such that $\rho$ is continuous in $\Omega$ and the space of all polynomials is dense
in $\cL^2(\Omega,\mu)$. Let $\ww$ be a $d$-tuple of positive real numbers.
We say that the triple $(\Omega,\LL,\mu)$
is a solution of the {\it Diffusion Orthogonal Polynomial Problem with weights $\ww$}
($\ww$-DOP problem for short) if, for any $n$, the space
$\{ P\in\bR[x_1,\dots,x_d]\mid \deg_{\ww}(P)\le n\}$ (considered as a subspace of
$\cL^2(\Omega,\mu)$) has an orthogonal basis formed by eigenvectors of $\LL$.
If in addition the equality (\eqSym) holds for any two smooth functions compactly
supported in $\bR^d$, then we say that $(\Omega,\LL,\mu)$ is a solution of
the {\it Strong} $\ww$-DOP problem ($\ww$-SDOP problem).
\medskip

 It is shown in
[\refBOZ, Prop.~2.11] that if $(\Omega,\LL,\mu)$ is a
solution of the $\ww$-SDOP problem, then $\LL$ is determined by the measure
density $\rho$ and the cometric $\g=(g^{ij})$ (the $g^{ij}$ are
the coefficients in (\eqLL)), namely,
$$
    \LL(f) = \frac1\rho\sum_{i,j}\partial_i
                     \Big( g^{ij}\rho\,\partial_jf\Big).            \eqno(\eqLLrev)
$$
We shall then speak also
of the triple $(\Omega,\g,\rho)$ as a solution of the $\ww$-SDOP problem.
Notice that (\eqLLrev) is the Laplace-Beltrami operator for the riemannian
metric $(g_{ij})=g^{-1}$ in the case when $\rho=(\det\g)^{-1/2}$.

As we mentioned above, any solution of the $\ww$-DOP problem with a bounded domain
$\Omega$ is a solution of the $\ww$-SDOP problem (see [\refBOZ, Prop.~\BOZpropSDOP]).
It is also shown in [\refBOZ] that any
solution of SDOP-problem is a solution of a certain algebraic problem formulated
in terms of the metric $g^{ij}$ only
(AlgDOP-problem; see~\S\sectGen\ below). This fact is proven in [\refBOZ] for the
usual degree but all the proofs extend without changes to any weighted degree as well.
Then, to find all two-dimensional solutions of the weighted SDOP problem,
we follow the same strategy as in [\refBOZ]. Namely, we first find
solutions of the algebraic problem over $\bC$
using some basic % prostejshie
properties of plane complex algebraic curves (see \S\S\sectBigW--\sectSmallW),
and then we look for $\Omega$ and $\rho$ (see \S\sectSDOP).

All the bounded domains $\Omega$ admitting a solution,
appeared already in the literature except, maybe, one
infinite family: the case (\modelBiv) with $m\ne n$ in Theorem~\thCompInf.
However orthogonal polynomials are being studied so long time in so many aspects
and the literature is so vast that we cannot be sure that this family is really new.

%%%%%%%%%%%%%%%%%%%%%%%%%%%%%%%%%%%%%%%%%%%%%%%%%%%%%%%%%%%%%%%%%
%%%%%%%%%%%%%%%%%%%%%%%%%%%%%%%%%%%%%%%%%%%%%%%%%%%%%%%%%%%%%%%%%
%%%%%%%%%%%%%%%%%%%%%%%%%%%%%%%%%%%%%%%%%%%%%%%%%%%%%%%%%%%%%%%%%
%%%%%%%%%%%%%%%%%%%%%%%%%%%%%%%%%%%%%%%%%%%%%%%%%%%%%%%%%%%%%%%%%
%%%%%%%%%%%%%%%%%%%%%%%%%%%%%%%%%%%%%%%%%%%%%%%%%%%%%%%%%%%%%%%%%

\head \sectGen. Some general facts about solutions of weighted
DOP/SDOP problem in an arbitrary dimension \endhead
\subhead\sectAlgDOP. The AlgDOP Problem
\endsubhead

Let $\ww=(w_1,\dots,w_d)$ and let
$(\Omega,\g,\rho)$ be a solution of the
$\ww$-SDOP problem in $\bR^d$. Let $\Delta=\det(g^{ij})$.
Let $\cP_{\ww}(n;\bK)$ be the vector space of polynomials
in $x_1,\dots,x_d$ with coefficients in $\bK$ whose $\ww$-weighted degree
is at most $n$. When $\bK$ is $\bR$, we write just $\cP_{\ww}(n)$.

Let $I(\partial\Omega)$ be the ideal in $\bR[x_1,\dots,x_d]$ 
of polynomials vanishing on $\partial\Omega$. The condition that
$\Omega$ is a natural domain implies that $I(\partial\Omega)$ is a principal ideal
(because $U\cap\partial\Omega$ cannot be of codimension $\ge2$ for any open $U$).
Let $\Gamma$ be a generator of $I(\partial\Omega)$, that is $\Gamma$
is a minimal polynomial vanishing on the boundary of $\Omega$.
In particular, if $\Gamma$ is not identically zero, then it is
{\it square-free}, i.e., it does not have multiple factors.
By convention we set $I(\varnothing)=\bR[x]$, i.e., $\Gamma=1$ in the case when
$\Omega$ is the whole $\bR^d$.

\proclaim{ Proposition \propAlgDOP } {\rm(See [\refBOZ, Thm.~2.21].)}
\roster
\item"(\condADi)"
    $g^{ij}\in\cP_{\ww}(w_i+w_j)$
    for any $i,j=1,\dots,d$.
    Hence $\Delta\in\cP_{\ww}(2w_1+\dots+2w_d)$.

\item"(\condADii)"
 $\partial\Omega\subset\{\Delta=0\}$, hence $\Gamma$ divides $\Delta$.

\item"(\condADiii)" 
For each $i=1,\dots,d$, one has
$$
      \sum_j g^{ij}\partial_j\Gamma = \Gamma S^i,
      \qquad S^i\in\cP_{\ww}(w_i).                        \eqno(\eqAlgDOP)
$$
\endroster
\endproclaim

This leads us to the following definition (cf.~[\refBOZ, Definition 3.2]).

\medskip\noindent
{\bf Definition \defAlgDOP. }
Let $\bK$ be $\bR$ or $\bC$ and let $\ww=(w_1,\dots,w_d)$
be a $d$-tuple of positive real numbers. A solution to the
{\it Algebraic Counterpart of the $\ww$-SDOP Problem over $\bK$}
($\ww$-AlgDOP Problem over $\bK$ for short)
is a pair $(\g,\Gamma)$ where $\g=(g^{ij})$ is a symmetric $d\times d$ matrix
with polynomial entries
and $\Gamma$ a polynomial such that
\roster
\item"(\condADi)"
     $g^{ij}\in\cP_{\ww}(w_i+w_j;\bK)$ for each $i,j=1,\dots,d$;
\item"(\condADii)"
     $\det\g$ is not identically zero,
      and $\Gamma$ is a square-free factor of $\det\g$;
\item"(\condADiii)"
     $\Gamma$ divides $\sum_i g^{ij}\partial_j\Gamma$
     for each $i=1,\dots,d$.
\endroster

\medskip
Thus Proposition~\propAlgDOP\ implies that if $(\Omega,\g,\rho)$
is a solution of the $\ww$-SDOP Problem and $\Gamma$
is a generator of $I(\partial\Omega)$, then $(\g,\Gamma)$ is a solution of
the $\ww$-AlgDOP Problem over $\bR$ and hence over $\bC$.

The following facts immediately follow
from the definition.

\proclaim{ Proposition \propW }
Let $\ww=(w_1,\dots,w_d)$ and $\ww'=(w'_1,\dots,w'_d)$ be two
$d$-tuples of positive weights. Suppose that
$(\g,\Gamma)$ is a solution of the $\ww$-AlgDOP Problem over $\bK$, and
$g^{ij}\in\cP_{\ww'}(w'_i+w'_j;\bK)$ for
each $i,j=1,\dots,d$. Then $(\g,\Gamma)$
is a solution of the $\ww'$-AlgDOP Problem over $\bK$.
\endproclaim

\proclaim{ Proposition \propGammaOne }
If $(\g,\Gamma)$ is a solution of the $\ww$-AlgDOP Problem over $\bK$
and $\Gamma_1$ is a factor of $\Gamma$, then $(\g,\Gamma_1)$
is also a solution of the $\ww$-AlgDOP Problem over $\bK$.
\endproclaim

If a square-free polynomial $\Gamma$ is given,
then it is easy to find all cometrics $\g$ such that
$(\g,\Gamma)$ is a solution of the $\ww$-AlgDOP problem.
Indeed, Condition (\condADiii) of Definition~\defAlgDOP\
(in the form (\eqAlgDOP)) provides a system
of linear equations for the coefficients of the polynomials $g^{ij}$ and $S^i$.
In \S\sectALGtoDOP\ we show how to find all $\rho$ for given $\Omega$ and $\g$.

%%%%%%%%%%%%%%%%%%%%%%%%%%%%%%%%%%%%%%%%%%%%%%%%%%%%%%%%%
%%%%%%%%%%%%%%%%%%%%%%%%%%%%%%%%%%%%%%%%%%%%%%%%%%%%%%%%%

\subhead\sectChangeVar. Admissible changes of variables
\endsubhead

A {\it $\ww$-admissible change of variables} is
a bijective polynomial mapping $\Phi:\bR^d\to\bR^d$,
$x\mapsto(y_1(x),\dots,y_d(x))$ with $\deg_{\ww}(y_i)=w_i$, $i=1,\dots,d$.
If $d=2$, then $(1,w)$-admissible automorphisms of $\bR^2$ for $w=1$ are
affine linear transformations, and for $w>1$, mappings of
the form
$$
   (x,y)\mapsto (\alpha x+\beta,\,\gamma y+p(x)),
             \qquad \alpha\gamma\ne0, \;\deg(p)\le w.        \eqno(\eqChangeVar)
$$
The following proposition is similar to [\refBOZ, Prop.~\BOZpropChangeVar]
and we omit its proof.

\proclaim{ Proposition \propChangeVar }
(a).
Let $(\Omega,\LL,\mu)$ be a solution of $\ww$-DOP (resp. $\ww$-SDOP) problem
and $\Phi$ be a $\ww$-admissible change of variables.
Let $\Omega_1=\Phi(\Omega)$, $\LL_1(f)=\LL(f\circ\Phi)\circ\Phi^{-1}$, and
$\mu_1(E)=\mu(\Phi^{-1}(E))$. Then $(\Omega_1,\LL_1,\mu_1)$ is also
a solution of $\ww$-DOP (resp. $\ww$-SDOP) problem.

\smallskip
(b). Let $(\g,\Gamma)$ be a solution of the $\ww$-AlgDOP Problem over $\bK$
and let $\Phi$ be a $\ww$-admissible change of variables.
Then $(\Phi_*(\g), \Gamma\circ\Phi^{-1})$ is also a solution
of the $\ww$-AlgDOP Problem over $\bK$.
\endproclaim

\medskip
\noindent
{\bf Example \exaChangeVar. }
Let $d=2$. For any $\ww=(w_1,w_2)$, the mapping
$\Phi:(x,y)\mapsto(x,-y)$ is $\ww$-admissible.
Thus, if $(\g,\Gamma)$ with
$\g=\left(\smallmatrix a&b\\b&c\endsmallmatrix\right)$
is a solution of the $\ww$-AlgDOP Problem over $\bK$, then
$(\Phi_*(\g),\Gamma(x,-y))$ with
$$
     \Phi_*(\g) = \left(\matrix a(x,-y) & -b(x,-y)\\ -b(x,-y) & c(x,-y)
          \endmatrix\right)
$$
is also a solution of the $\ww$-AlgDOP problem over $\bK$.

\medskip
\noindent
{\bf Example \exaChangeVarII. }
More generally, let still $d=2$. Any $(1,w)$-admissible
change of variables for $w>1$ is of the form (\eqChangeVar) and it
is a composition
of the following variable changes $(x,y)\mapsto(X,Y)$:
$$
   T:(x,y)\mapsto(x+\beta,y),\quad
   H:(x,y)\mapsto(\alpha x,\gamma y), \quad
   S:(x,y)\mapsto(x,y+p(x)).
$$
They transform $\g=\left(\smallmatrix a&b\\b&c\endsmallmatrix\right)$ as follows
(see [\refBOZ, eq.~(\BOZeqChangeVar)]): $T_*(\g) = \g(X-\beta,Y-\beta)$,
$$
   H_*(\g) = \left(\matrix
    a\alpha^2     &  b\alpha\gamma \\
    b\alpha\gamma &  c\gamma^2\endmatrix\right)
    _{\smallmatrix x=X/\alpha\\y=Y/\gamma\endsmallmatrix}\,
\quad
   S_*(\g) = \left(\matrix
     a       &  p'a + b\\
     p'a + b &  (p')^2a + 2p'b + c\endmatrix\right)
    _{\smallmatrix\hskip-22pt x=X\\y=Y-p(X)\endsmallmatrix}
$$

%%%%%%%%%%%%%%%%%%%%%%%%%%%%%%%%%%%%%%%%%%%%%%%%%%%%%%%%%
%%%%%%%%%%%%%%%%%%%%%%%%%%%%%%%%%%%%%%%%%%%%%%%%%%%%%%%%%

\subhead\sectALGtoDOP. From AlgDOP to DOP/SDOP
\endsubhead

As we told in \S\sectAlgDOP, any solution $(\Omega,\g,\mu)$ of the $\ww$-SDOP
problem gives a solution $(\Gamma,\g)$ of the $\ww$-AlgDOP problem where
$\Gamma$ is a generator of the ideal $I(\partial\Omega)$.
In this subsection we discuss how to find all possible $(\Omega,\g,\mu)$
from a given $(\Gamma,\g)$.

So, let $(\Gamma,\g)$ be a solution of the $\ww$-AlgDOP problem over $\bR$.
First of all we find all connected components $\Omega$ of $\bR^d\setminus\{\Gamma=0\}$
such that $\g$ is positive definite on $\Omega$. Then it remains to find all
measure densities $\rho$
such that the operator $\LL$ given by (\eqLLrev) is of the form (\eqLL)
with $b^i\in\cP_\ww(w_i)$ for each $i=1,\dots,d$. By comparing (\eqLL) with (\eqLLrev)
we obtain
$$
      b^i = \sum_j \partial_j g^{ij} + \sum_j g^{ij}\partial_j h, \qquad h=\log\rho.
$$
Hence
$$
      \partial_j h = \sum_i g_{ij} L^i,  \qquad
      L^i = b^i - \sum_j \partial_j g^{ij},                        \eqno(\eqDh)
$$
where $(g_{ij})=\g^{-1}$. We have
$b^i\in\cP_\ww(w_i)\Leftrightarrow L^i\in\cP_\ww(w_i)$. The identities
$\partial_k(\partial_j h)=\partial_j(\partial_k h)$ combined with (\eqDh) yield
$$
      \partial_k\Big(\sum_i g_{ij} L^i\Big)
     =\partial_j\Big(\sum_i g_{ik} L^i\Big)                        \eqno(\eqDkDj)
$$
which is a system of linear equations for the coefficients of all $L^i$.
These observations lead to the following algorithm of finding all solutions for $\rho$.
We start with polynomials $L^i$, $\deg_\ww L^i=w_i$, whose coefficients we compute
by solving the system of linear equations coming from (\eqDkDj).
The solution may depend on several parameters.
Then we find $h$ by integrating the $L^i$'s and set $\rho=\exp(h)$.
Finally we choose the values of the parameters such that $Q(x)\rho(x)$,
is integrable over $\Omega$ for any polynomial $Q$
(when $\Omega$ is bounded, it is enough to demand that $\int_\Omega\rho\,dx<\infty$).

These observations allow us to prove the following fact.

\proclaim{ Proposition \propOPlus }
Let $(\Omega,\g,\rho)$ be a solution of the $\ww$-SDOP problem in $\bR^d$
such that $\g=\diag\big(g^{11}(x_1),\dots,g^{dd}(x_d)\big)$. Then there is a partition
of the set of indices $\{1,\dots,d\}=J_1\sqcup\dots\sqcup J_m$,
$J_k=\{j_k(1),\dots,j_k(d_k)\}$, such that:
\roster
\item"$\bullet$" $\rho=\rho_1(\xx_1)\dots\rho_m(\xx_m)$,
                 where $\xx_k=(x_{j_k(1)},\dots,x_{j_k(d_k)})$,
                 $k=1,\dots,m$;
\item"$\bullet$" if $d_k>1$, then $w_i=w_j$ and $g^{ii}$ is constant
                 for any $i,j\in J_k$.
\endroster 
\endproclaim

\demo{ Proof } We consider only the case $d=2$.
It is not difficult to derive from it the general case. Let $h=\log\rho$.
It is enough to prove that either $\g$ is constant, or $\partial_{12}h=0$ and then
$h=h_1(x_1)+h_2(x_2)$ whence $\rho=\rho_1(x_1)\rho_2(x_2)$.
Indeed, the inverse matrix of $g$ is $\diag(g_{11},g_{22})$ where $g_{ii}=1/g^{ii}$.
Then by (\eqDkDj) we have
$$
     \partial_{12}h = \frac{\partial_2 L^1(x_1,x_2)}{g^{11}(x_1)}
                    = \frac{\partial_1 L^2(x_1,x_2)}{g^{22}(x_2)},
           \qquad \deg_\ww L^i\le w_i.
$$
Therefore $\partial_{12}h$ is a polynomial, hence
$g^{ii}$ divides $\partial_j L^i$ for $i\ne j$.
Thus, if $\partial_{12}h$ is nonzero, we have
$0\le\deg_\ww g^{11}\le\deg_\ww\partial_2 L^1=\deg_\ww L^1-w_2\le w_1 - w_2$
and, symmetrically, $0\le\deg_\ww g^{22}\le w_2-w_1$, which implies
$w_1=w_2$ and $\deg_\ww\g=0$, i.e., $\g$ is constant.
\qed
\enddemo

To conclude this subsection, we give some conditions (necessary or sufficient)
on $\rho$ to be compatible with given
$\g$ and $\Omega$. They are proven in [\refBOZ] for the usual degree
but the proofs can be easily adapted for any weighted degree.

\proclaim{ Proposition \propDOP } {\rm(See [\refBOZ, Thm.~2.21].)}
Let $(\g,\Gamma)$, $\Gamma=\Gamma_1\dots\Gamma_s$,
be a solution of $\ww$-AlgDOP problem over $\bR$
(recall that then $\Gamma$ is squarefree).
Let $\Omega$ be a bounded domain such that $\partial\Omega\subset\{\Gamma=0\}$
and $\g$ is positive definite in $\Omega$.
Let $a_1,\dots,a_s$ be real numbers such that $\mu(\Omega)<\infty$ where
$\mu$ is the measure with density $\rho=\prod_\nu\Gamma_\nu^{a_\nu}$
(for example, $a_\nu\ge 0$ for each $\nu$). Then $(\Omega,\LL,\mu)$ is a solution of
$\ww$-DOP problem, where $\LL$ is given by (\eqLLrev).
\endproclaim

\noindent{\bf Definition \defMaxBndry.}
A solution $(\g,\Gamma)$ of the $\ww$-AlgDOP problem is called
{\it maximal} (and in this case $\Gamma$ is called
a {\it maximal boundary for} $\g$)
if $\Gamma_1$ divides $\Gamma$ for any solution $(\g,\Gamma_1)$.
By Proposition~\propGammaOne, a maximal solution
is unique for any given $\g$.

\proclaim{ Proposition \propMes } {\rm(See [\refBOZ, Props.~2.15, 2.17].)}
Let $(\Omega,\g,\rho)$ be a solution of $\ww$-SDOP problem in $\bR^d$
and $\Delta=\det(\g)$. Let $\Gamma$ be the maximal boundary for $\g$ and let
$\Gamma_1,\dots,\Gamma_s$ be its irreducible (over $\bC$) factors.
Suppose that each factor $\Gamma_k$
occurs in $\Delta$ with multiplicity $1$, i.e.,
$\Gamma_k^2$ does not divide $\Delta$.
Then (see Remark~\remMes)
$$
   \rho=\Gamma_1^{p_1}\dots\Gamma_s^{p_s}\exp(Q)           \eqno(\eqMes)
$$
for some $p_1,\dots,p_s\in\bC$ and a polynomial $Q$ such that, for each $j=1,\dots,d$,
$$
   w_j\deg_{x_j}(Q\Delta) \le 2w_1+\dots+2w_d.              \eqno(\eqMesIneq)
$$
\endproclaim

\noindent{\bf Remark \remMes.} In Proposition~\propMes,
if $\lambda\Gamma_k$, $\lambda\in\bC$, is real for some $k$,
we always assume that $\Gamma_k$ is real
and positive on $\Omega$; in this case $p_k$ is real. 
 Otherwise $\bar\Gamma_k$ is also a factor of $\Gamma$
 and it must occur in $\rho$ with the power $\bar p_k$ because $\rho$ is real.
 In this case $\Gamma_k^{p_k}\bar\Gamma_k^{\bar p_k}$
 is understood as a single-valued branch of this function on $\Omega$.
 Notice that choosing another single-valued branch
 we just change the constant term of $Q$.

\medskip\noindent
{\bf Remark \remMesBis.}
  Proposition~\propMes\ admits the following refinement.
  We still have (\eqMes) and (\eqMesIneq) for a variable $x_j$ even
  when there are multiple factor $\Gamma_k$ of $\Delta$ but they
  are polynomials in variables $(x_i)_{i\in I}$ not including $x_j$,
  that is $j\not\in I$.
  In this case $Q$ is a polynomial in $(x_i)_{i\not\in I}$ whose
  coefficients are rational functions in $(x_i)_{i\in I}$.

\proclaim{ Corollary \corMes }
Under the hypothesis of Proposition~\propMes, assume that $mw_d=2(w_1+\dots+w_d)$
and $w_d=w_i n_i$, $i=1,\dots,d$, for some $m,\,n_1,\dots,n_d$
{\rm(e.g.~$d=2$ and $(w_1,w_2)=(1,2)$).}
Suppose that $\deg_{x_d}\Delta=m$. Then $Q=\const$ in (\eqMes).
\endproclaim

\demo{ Proof } By an admissible change of variables
$y_d=x_d+\sum_{i=1}^{d-1}a_i x_i^{n_i}$ and $y_j=x_j$ for $j<d$,
we may achieve that $w_j\deg_{y_j}\Delta=2\sum w_i$ for each $j$.
\qed\enddemo

Recall that $(g_{ij})=\g^{-1}$, hence $g_{ij}=\hat g_{ij}/\det g$ with
a polynomial $\hat g_{ij}$.

\proclaim{ Proposition \propMesII } {\rm(See [\refBOZ, Cor.~2.19].)}
Let $U=\bR\times\bB^{d-1}=\{x_2^2+\dots+x_d^2<1\}$ and
$U_+=\bR_+\times\bB^{d-1}=U\cap\{x_1>0\}$.
Let $M_1=\max_j\big(\lfloor w_j/w_1\rfloor + \deg_{x_1}\hat g_{1j}\big)$.
Let $(\Omega,\g,\rho)$ be a solution
of the $\ww$-SDOP problem, $\Delta=\det\g$. Suppose
that $U \subset\Omega$ (resp. $U_+ \subset\Omega$).
Then $\deg_{x_1}\Delta<M_1$ (resp. $\deg_{x_1}\Delta<1+M_1$).
\endproclaim

%%%%%%%%%%%%%%%%%%%%%%%%%%%%%%%%%%%%%%%%%%%%%%%%%%%%%%%%%
%%%%%%%%%%%%%%%%%%%%%%%%%%%%%%%%%%%%%%%%%%%%%%%%%%%%%%%%%
%%%%%%%%%%%%%%%%%%%%%%%%%%%%%%%%%%%%%%%%%%%%%%%%%%%%%%%%%
%%%%%%%%%%%%%%%%%%%%%%%%%%%%%%%%%%%%%%%%%%%%%%%%%%%%%%%%%
%%%%%%%%%%%%%%%%%%%%%%%%%%%%%%%%%%%%%%%%%%%%%%%%%%%%%%%%%
%%%%%%%%%%%%%%%%%%%%%%%%%%%%%%%%%%%%%%%%%%%%%%%%%%%%%%%%%
%%%%%%%%%%%%%%%%%%%%%%%%%%%%%%%%%%%%%%%%%%%%%%%%%%%%%%%%%
%%%%%%%%%%%%%%%%%%%%%%%%%%%%%%%%%%%%%%%%%%%%%%%%%%%%%%%%%
%%%%%%%%%%%%%%%%%%%%%%%%%%%%%%%%%%%%%%%%%%%%%%%%%%%%%%%%%
%%%%%%%%%%%%%%%%%%%%%%%%%%%%%%%%%%%%%%%%%%%%%%%%%%%%%%%%%

\head \sectDimTwoAlg. Weighted AlgDOP Problem in $\bC^2$
\endhead

\subhead\sectDimTwoAlg.1. Notation
\endsubhead
In \S\S\sectDimTwoAlg--\sectSmallW\
we study the $\ww$-AlgDOP Problem over $\bC$
in dimension 2 for any pair of weights $\ww$.
It is clear that a multiplication of $\ww$ by a positive number does not change
the problem. Therefore we assume throughout this section that
$(\g,\Gamma)$, is a solution of the $\ww$-AlgDOP Problem over $\bC$ for $\ww=(1,w)$
with a real $w>1$ (the case $w=1$ is already done in [\refBOZ]).

We denote variables by $(x,y)$ and we set
$\g=\left(\smallmatrix a&b\\b&c\endsmallmatrix\right)$.
We denote the coefficient of $y^m$ in $a(x,y)$ by $a_m(x)$ and
the coefficient of $x^ky^m$ by $a_{km}$ and we use similar notation for $b$ and $c$.
Sometimes we write $(a,b,c;\Gamma)$ instead of $(\g,\Gamma)$ when
speaking of a solution of the $\ww$-AlgDOP Problem.

As usual, given a polynomial $P(x,y)=\sum p_{km} x^k y^m$, we define its
Newton polygon $\cN(P)$ as the convex hull in $\bR^2$ of the finite set
$\{(k,m)\mid a_{km}\ne 0\}$.
Condition (\condADi) of Definition~\defAlgDOP\ means that the Newton polygons
of $a$, $b$, $c$, and $\Delta$ are contained in the polygons shown in
Figures~\figABCi--\figABCii.

\midinsert
\centerline{\epsfxsize=80mm\epsfbox{abc1.eps}}
\centerline{$a$ \hskip 12mm $b$ \hskip 16mm $c$ \hskip 20mm $\Delta$ \hskip 10mm}
\botcaption{Figure~\figABCi}
      Polygons containing $\cN(a)$, $\cN(b)$, $\cN(c)$, $\cN(\Delta)$
       for $1<w\le2$.
\endcaption
\endinsert

\midinsert
\centerline{\epsfxsize=90mm\epsfbox{abc2.eps}}
\vskip -20pt
\centerline{\hskip 31mm ${}_{1+w}$ \hskip 21mm ${}_{2w}$
             \hskip 28mm ${}_{2+2w}$}
\vskip 8pt
\centerline{$a$ \hskip 12mm $b$ \hskip 18mm $c$
          \hskip 24mm $\Delta$ \hskip 15mm}
\botcaption{Figure~\figABCii}
      Polygons containing $\cN(a)$, $\cN(b)$, $\cN(c)$, $\cN(\Delta)$
       for $w>2$.
\endcaption
\endinsert

%%%%%%%%%%%%%%%%%%%%%%%%%%%%%%%%%%%%%%%%%%%%%%%%%%%%%%%%%
%%%%%%%%%%%%%%%%%%%%%%%%%%%%%%%%%%%%%%%%%%%%%%%%%%%%%%%%%

\medskip
\subhead \sectW. Change of the weights and the $(1,\infty)$-AlgDOP Problem
\endsubhead

According to Proposition~\propW, any solution of $(1,w)$-AlgDOP Problem
with $1<w\le 2$ is a solution of the corresponding $(1,2)$-AlgDOP Problem
(see Figure~\figABCi).
Similarly (see Figure~\figABCii) any solution of $(1,w)$-AlgDOP Problem with $w>2$
is a solution of the corresponding $(1,w')$-Problem for any $w'>w$.

So, we say that $(\g,\Gamma)$ is a
solution of the $(1,\infty)$-AlgDOP Problem, if it is a solution of the
$(1,w)$-AlgDOP Problem
for some $w>2$, and a coordinate change (\eqChangeVar) with an arbitrary
polynomial $p(x)$ will be called $(1,\infty)$-admissible.

In the next two sections we find all solutions of the
$(1,\infty)$- and $(1,2)$-AlgDOP Problems up to $(1,\infty)$- and $(1,2)$-admissible
coordinate change. 

%%%%%%%%%%%%%%%%%%%%%%%%%%%%%%%%%%%%%%%%%%%%%%%%%%%%%%%%%
%%%%%%%%%%%%%%%%%%%%%%%%%%%%%%%%%%%%%%%%%%%%%%%%%%%%%%%%%

\medskip
\subhead \sectLocal. Local branches of the curve $\Gamma=0$
\endsubhead

Let $P(x,y)$ be a squarefree polynomial. A {\it local branch} of $P$
(or, equivalently, of the curve $P=0$)
is a pair $\gamma=(\xi,\eta)$ of germs at $0$ of meromorphic functions
such that $P(\xi(t),\eta(t))$ is identically zero.
Any meromorphic germ $t\mapsto\gamma(t)=(\xi(t),\eta(t))$ defines
a {\it valuation} $v_\gamma:\bC[x,y]\to\bZ\cup\{\infty\}$,
$v_\gamma(Q)=\ord_t Q(\xi(t),\eta(t))$ where $\ord_t(0)=\infty$ and
$\ord_t f(t)=p_m$ if $f(t)=\sum_{k\ge m}p_kt^k$ and $p_m\ne 0$.

Let $(a,b,c;\,\Gamma)$ be a solution of the $\ww$-AlgDOP problem over $\bC$ for
some $\ww=(1,w)$. Condition (\condADiii) of Definition~\defAlgDOP\ reads
$$
   a\diff\Gamma{x} + b\diff\Gamma{y} = L_1\Gamma, \qquad
   b\diff\Gamma{x} + c\diff\Gamma{y} = L_2\Gamma.         \eqno(\eqBord)
$$
It is easy to check that this condition is equivalent to
$$
   b(\xi,\eta)\dot\xi = a(\xi,\eta)\dot\eta, \qquad
   c(\xi,\eta)\dot\xi = b(\xi,\eta)\dot\eta
                                                                \eqno(\eqParam)
$$
for any local branch of $\Gamma$. Condition (\eqParam) implies that
$$
  v_\gamma(a) - v_\gamma(b) = v_\gamma(b) - v_\gamma(c)
      = \ord_t(\dot\xi) - \ord_t(\dot\eta)                 \eqno(\eqValu)
$$
if both $\xi(t)$ and $\eta(t)$ are non-constant (see [\refBOZ, Lemma~3.3]).

The following fact is well-known and immediately follows from the definitions.

\proclaim{ Lemma \lemNewton } Let $F$ be a polynomial in $(x,y)$ and
let $(p,q)\in\bZ^2\setminus\{(0,0)\}$.
A local branch $\gamma$ of $F$ such that $\ord_t(\gamma)=(p,q)$ exists
if and only if the vector $(p,q)$ is orthogonal to some edge of $\cN(F)$
and points from this edge inward $\cN(\Gamma)$.
\qed\endproclaim

 \proclaim{ Lemma \lemChangeW }
 Let $(a,b,c;\,\Gamma)$ is a solution of $(1,w)$-AlgDOP Problem with $w>1$.
 Suppose that
 $\Gamma$ is divisible by $x$. Then $a$ and $b$ are divisible by $x$, hence
 $\deg_y a\le 1$ and $(a,b,c;\,\Gamma)$ is a solution of the
 $(1,\infty)$-AlgDOP Problem.
 \endproclaim

 \demo{ Proof }
 We have $a(0,t)=b(0,t)=0$ by (\eqParam), hence $x$ is a factor of $a$ and $b$
 whence $\deg_y a\le 1$ (see Figure~\figABCi).
 \qed\enddemo

\proclaim{ Lemma \lemValuA }
Let $w>1$. Suppose that $\Gamma$ has a branch
$\gamma=(\xi(t),\eta(t))$ such that
$v_\gamma(x,y)=\ord_t\gamma=(p,q)$
with
$q<0<p$. Then $b_1 = b_{11} x$ and $a=a_{20} x^2$, hence
$(\g;\,\Gamma)$ is a solution of $(1,\infty)$-AlgDOP Problem.
\endproclaim

\demo{Proof}
We have $\ord_t\dot\xi-\ord_t\dot\eta=p-q$, hence
$v_\gamma(a)-v_\gamma(b) = v_\gamma(b) - v_\gamma(c) = p-q$ by (\eqValu).
On the other hand we have $v_\gamma(c)\ge v_\gamma(y^2)=2q$, hence
$$
     v_\gamma(b) = v_\gamma(c)+p-q\ge p+q.       \eqno(\eqPfLemA)
$$
If $b_{01}$ were non-zero, then we would have
$v_\gamma(b)=v_\gamma(y)=q$ which contradicts (\eqPfLemA) because $p>0$.
Hence $b_1=b_{11}x$ (as required), and this fact implies that
$$
   v_\gamma(b)\ge\min(v_\gamma(1),v_\gamma(xy))=\min(0,p+q),
$$
thus $v_\gamma(a)=v_\gamma(b)+p-q\ge\min(0,p+q)+p-q= p+\min(-q,p)>p$.
Since the only monomials which may occur in $a$ are $(y,1,x,x^2)$
and their $\gamma$-valuations are $(q,0,p,2p)$,
we see that $a=a_{20}x^2$.
\qed\enddemo

\proclaim{Lemma \lemValu } Let $1<w\le 2$.
Let $\gamma=(\xi(t),\eta(t)$ be a local branch of $\Gamma$, and
$(p,q)=\ord_t\gamma=v_\gamma(x,y)$.

\smallskip
(a). If $(p,q)=(2,1)$ then $a_{00}=a_{01}=b_{00}=0$,
     and hence $\deg_y\Delta\le2$ and
     $(a,b,c;\Gamma)$ is a solution of the $(1,\infty)$-AlgDOP problem.

\smallskip
(b). If $q>0$, then $p\ne 3$.

\smallskip
(c). If $p=-3$ and $-4\le q\le -1$, then $\deg b_0\le2$, $\deg c_0\le 2$,
     and $\deg c_1\le 1$, hence
     $(a,b,c;\Gamma)$ is a solution of the $(1,1)$-AlgDOP problem.

\smallskip
(d). If $q<2p<0$, then $a_{01}=0$, hence
     $(a,b,c;\Gamma)$ is a solution of the $(1,\infty)$-AlgDOP problem.
\endproclaim

\demo{ Proof }
(a). If $(p,q)=(2,1)$, then $\ord_t(\dot\xi,\dot\eta)=(1,0)$, hence
by (\eqValu) we have $v_\gamma(b)=v_\gamma(c)+1\ge 1$ and
$v_\gamma(a)=v_\gamma(b)+1=v_\gamma(c)+2\ge 2$, and the result follows from the fact
that $1$ (resp. $y$) is the only monomial of $a$ and $b$ of the
$\gamma$-valuation equal to $0$ (resp. to $1$).
The condition $a_{01}=0$ implies that 
$(a,b,c;\Gamma)$ is a solution of the $(1,1)$-AlgDOP problem (see Figure~\figABCii).

\smallskip
(b). Suppose that $p=3$ and $q>0$,
Then by the arguments similar to those in (a), we are going to show that
$\deg_y a\le0$ and hence
$\deg_y\Delta\le2$ which contradicts the condition $\ord_t\eta=3$.
Notice that $q\le\deg_x\Delta\le 6$.

If $q=1$, the proof is the same as in (a).

If $q=2$, then
by (\eqValu) we have
$v_\gamma(a)=v_\gamma(b)+1$.
Since $v_\gamma(a)\ne 1$ and $v_\gamma(b)\ne 1$, this
is possible only when $v_\gamma(a)\ge 3$. By combining this fact with
$v_\gamma(1,y,x,x^2)=(0,2,3,6)$, we obtain
$a_{00}=a_{01}=0$ whence $\deg_y a\le 0$.

The case $q=3$ is reduced to $q>3$ by a change of variables
$(x,y)\mapsto(x,y+\lambda x)$.

If $q=4$, then
$v_\gamma(a)\in\{0,3,4,6\}$,
$v_\gamma(b)\in\{0,3,4,6,7,9\}$,
$v_\gamma(c)\in\{0,3,4,6,\dots\}$.
Then one can check that $v_\gamma(c)-v_\gamma(b)=v_\gamma(b)-v_\gamma(a)=1$
(this is (\eqValu))
is possible only when $v_\gamma(a)=6$, hence $\deg_y a\le 0$.

If $q=5$, then
$v_\gamma(a)\in\{0,3,5,6\}$,
$v_\gamma(b)\in\{0,3,5,6,8,9\}$,
$v_\gamma(c)\in\{0,3,5,6,8,\dots\}$.
Then $v_\gamma(c)-v_\gamma(b)=v_\gamma(b)-v_\gamma(a)=2$
(this is (\eqValu))
is possible only when $v_\gamma(a)=6$, hence $\deg_y a\le 0$.

The case $q=6$ is reduced to $q>6$ by a change of variables
$(x,y)\mapsto(x,y+\lambda x^2)$.

\smallskip
(c). We need to show that
of $x^3$, $x^4$, and $x^2y$ do not occur in $b$ and in $c$.

If $p=-3$ and $-3\le q\le -1$, then (\eqValu) implies
$v_\gamma(c)\ge v_\gamma(b)\ge v_\gamma(a)\ge v_\gamma(x^2)=-6$.
On the other hand,
$v_\gamma(x^2y)=-6+q\ge -7$, $v_\gamma(x^3)=-9$, $v_\gamma(x^4)=-12$,
and these are the only monomials which might occur in $b$ and in $c$ with the
respective $\gamma$-valuations.

If $(p,q)=(-3,-4)$, then  (\eqValu) implies $v_\gamma(a)\ge v_\gamma(x^2)=-6$,
$v_\gamma(b)=v_\gamma(a)-1\ge-7$, and
$v_\gamma(c)=v_\gamma(b)-1\ge-8$.
On the other hand,
$v_\gamma(x^2y)=-10$, $v_\gamma(x^3)=-9$, $v_\gamma(x^4)=-12$,
and these are the only monomials which might occur in $b$ and in $c$ with the
respective $\gamma$-valuations.

\smallskip
(d). The condition $q<2p<0$ implies $v_\gamma(b)\ge v_\gamma(xy)= p+q$.
It implies also $p-q>-p$.
By (\eqValu) we have $v_\gamma(a)=v_\gamma(b)+p-q$. Hence
$$
  v_\gamma(a)=v_\gamma(b)+p-q>v_\gamma(b)-p\ge(p+q)-p=q=v_\gamma(y).
$$
Therefore $a_{01}=0$ because $y$ is the only monomial which may occur in $a$
with the $\gamma$-valuation equal to $q$.
The condition $a_{01}=0$ implies that 
$(a,b,c;\Gamma)$ is a solution of the $(1,1)$-AlgDOP problem (see Figure~\figABCii).
\qed\enddemo

%%%%%%%%%%%%%%%%%%%%%%%%%%%%%%%%%%%%%%%%%%%%%%%%%%%%%%%%%
%%%%%%%%%%%%%%%%%%%%%%%%%%%%%%%%%%%%%%%%%%%%%%%%%%%%%%%%%
%%%%%%%%%%%%%%%%%%%%%%%%%%%%%%%%%%%%%%%%%%%%%%%%%%%%%%%%%
%%%%%%%%%%%%%%%%%%%%%%%%%%%%%%%%%%%%%%%%%%%%%%%%%%%%%%%%%
%%%%%%%%%%%%%%%%%%%%%%%%%%%%%%%%%%%%%%%%%%%%%%%%%%%%%%%%%
%%%%%%%%%%%%%%%%%%%%%%%%%%%%%%%%%%%%%%%%%%%%%%%%%%%%%%%%%
%%%%%%%%%%%%%%%%%%%%%%%%%%%%%%%%%%%%%%%%%%%%%%%%%%%%%%%%%
%%%%%%%%%%%%%%%%%%%%%%%%%%%%%%%%%%%%%%%%%%%%%%%%%%%%%%%%%

\head\sectBigW. Solution of the $(1,w)$-AlgDOP Problem in $\bC^2$ for $w>2$
\endhead

In this section we give all solutions of the $(1,\infty)$-AlgDOP Problem over $\bC$
up to $(1,\infty)$-admissible change of coordinates.
As we explained in \S\sectW, this gives all solutions of
the $(1,w)$-AlgDOP Problem for all $w>2$.

\proclaim{ Lemma \lemA }
Let $(a,b,c;\Gamma)$ be a solution of the $(1,w)$-AlgDOP problem over $\bC$ with
$w>1$. Suppose that $\deg_y\Gamma=2$ and that $\Gamma$
is not divisible by any non-constant polynomial in $x$.
Then $\Gamma$ is monic with respect to $y$, i.e., $y^2$ is the only monomial
of $\Gamma$ whose $y$-degree is $2$.
\endproclaim

\demo{Proof} First, remark that $a\ne 0$. Indeed, otherwise $\Delta=b^2$, hence
$\Gamma$ would be a square-free factor of $b$ which contradicts the condition
$\deg_y\Gamma=2$.
Suppose that the coefficient of $y^2$ in $\Gamma$ is a non-constant polynomial in $x$.
Without loss of generality we may assume that one of its roots is $0$.
Then $\Gamma$ has a branch $\gamma=(\xi(t),\eta(t))$ such that
$\ord_t\eta<0<\ord_t\xi$. Then, by Lemma~{\lemValuA}, we have
$a=a_{20}x^2$ and $b=b_{11}xy+b_0(x)$.
Since $a\ne0$, we may assume that $a=x^2$.

If $b_{00}=0$, then $x^2$ is a factor of $ac$ and of $b^2$, hence $x^2$ is
a factor of $\Delta$. Since $\Delta$ may not have
have monomials $y^2x^k$ with $k>2$, it follows that $\Gamma$ is monic in $y$.

Now consider the case $b_{00}\ne 0$. Then the constant term of $\Delta$
is $-b_{00}^2\ne 0$. 
We have $a=x^2$ and $b_1=b_{11}x$, hence $x^2y^2$ is the only monomial of
$y$-degree $2$ which may occur in $\Delta$. Since $\deg_y\Delta=2$,
we conclude that $\Delta$ cannot be divisible by any non-constant polynomial in $x$,
i.e., $\Gamma=\Delta$.

The coefficients of $y^2$ in $\Delta$ and in
$a\diff\Delta{x}+b\diff\Delta{y}$
are $d_{22}x^2$ and $2d_{22}(1+b_{11})x^3$ respectively where
$d_{22}=c_{02}-b_{11}^2$.
Then (\eqBord) implies $L_1=2(1+b_{11})x$.
By plugging this back into (\eqBord), we obtain a contradiction. Indeed, we have:
$$
\split
   \Delta&=-b^2+O(x^2)
         =-b_{00}^2 - 2b_{00}(b_{10}+b_{11}y)x+O(x^2),
\\
   b\diff\Delta{y} &= \big(b_{00}+O(x)\big)
                         \big(-2b_{00}b_{11}x + O(x^2)\big)
                      = -2b_{00}^2b_{11}x+O(x^2),
\\
   L_1\Delta& = -2b_{00}^2(1+b_{11})x + O(x^2).
\endsplit
$$
Hence
$a\diff\Delta{x}+b\diff\Delta{y}-L_1\Delta = 2b_{00}^2 x+O(x^2)\ne 0$.
\qed\enddemo

\proclaim{ Lemma \lemSym }
Let $\Gamma=y^2-p(x)$ and
let $(a,b,c;\,\Gamma)$ be a solution of the $(1,w)$-AlgDOP problem with $w>2$.
Then $b_0=c_1=0$, i.e.,
$a$ and $c$ are even with respect to $y$,
and $b$ is odd with respect to $y$ (this means that
the corresponding metric is invariant under the symmetry $y\mapsto-y$).
\endproclaim

\demo{ Proof }
Let us set $\hat a=a$, $\hat b(x,y) = -b(x,-y)$, $\hat c(x,y)=c(x,-y)$,
and $\hat\g=(\hat a,\hat b,\hat c)$.
Since $\Gamma$ is symmetric, Proposition~{\propChangeVar} implies that
$(\hat\g;\,\Gamma)$ is also a solution to the
$(1,w)$-AlgDOP problem (see Example~\exaChangeVar).
Hence both $(\g;\Gamma)$ and $(\hat\g;\Gamma)$ satisfy 
the equations (\eqBord) and, by linearity,
$\big(\frac12(\g-\hat\g);\,\Gamma\big)$ satisfies it as well.
Since $\frac12(\g-\hat\g)=(0,b_0,c_1y)$, this means
$2yb_0 = (y^2+p(x))L_1$ and  $-b_0p'(x)+2y^2c_1 = (y^2+p(x))L_2$.
The first equation implies $b_0=0$. Plugging this into the second equation,
we obtain $2y^2c_1 = (y^2+p(x))L_2$.
Note that $p\ne0$ because $\Gamma$ cannot have multiple factors.
Hence the last equations implies $c_1=0$, i.e., $\g-\hat\g=0$.
\qed\enddemo

\proclaim{ Proposition \propDegGaTwo }
The following is a complete list of solutions $(\g;\Gamma)$ of
the $(1,\infty)$-AlgDOP problem over $\bC$ up to $(1,\infty)$-admissible change of
variables under condition that $\deg_y\Gamma=2$:
\roster
\item"(i)" $\Gamma=(1-x)^m(1+x)^n - y^2$, $m,n\ge 1$,
$$
  g=\left(\matrix 1-x^2 & \tfrac12\big((n-m)-(n+m)x\big)y \lower3pt\hbox{\mathstrut}\\
    * & \frac14\big((n-m)-(n+m)x\big)^2(1-x)^{m-1}(1+x)^{n-1} - c_{02}\Gamma
    \endmatrix\right)\,;
$$

\smallskip
\item"(ii)" $\Gamma=x^n-y^2$, $n\ge 1$,
$$
  g=\left(\matrix x & \;\tfrac12ny \lower5pt\hbox{\mathstrut}\\
        \tfrac12ny  & \;\frac14n^2 x^{n-1}-c_{02}\Gamma
    \endmatrix\right)\,;
$$

\item"(iii)" $\Gamma=\Gamma_2x^k$ where
             $\Gamma_2=(x_0-x)^n-y^2$, $n\ge 0;$ $k,x_0\in\{0,1\}$, $(n,c_{02})\ne(0,0)$,

$$
  g=\left(\matrix x(x_0-x) & -\tfrac12nxy \lower5pt\hbox{\mathstrut}\\
        -\tfrac12nxy       & \frac14n^{2}x(x_0-x)^{n-1} - c_{02}\Gamma_2
    \endmatrix\right)\,;
$$

\item"(iv)"$\Gamma=    x^k(1-y^2)$, $\g=\diag(  x^k,-c_{02}(1-y^2))$, $k\in\{0,1\}$,
                                                                      $c_{02}\ne 0$;
\item"(v)" $\Gamma=(1-x^2)(1-y^2)$, $\g=\diag(1-x^2,-c_{02}(1-y^2))$, $c_{02}\ne 0$.
\endroster
\endproclaim

\demo{ Proof }
Let $\Gamma=\Gamma_0\Gamma_2$ where $\deg_y\Gamma_0=0$,
$\deg_y\Gamma_2=2$, and $\Gamma_2$ is not divisible by any non-constant
polynomial in $x$. Then $\Gamma_2$ is monic in $y$ by Lemma~\lemA\ combined with
Proposition~\propGammaOne, hence
it can be reduced to the form $\Gamma_2=y^2+p(x)$ by a $(1,\infty)$-admissible
change of coordinates.
Since $\Gamma$ cannot have any multiple factor,
$p$ is not identically zero. Moreover, by rescaling the $y$-coordinate, we may
set the leading coefficient of $p$ to any prescribed nonzero complex number.
By Lemma~{\lemSym}, we have
$a=a(x)$, $b=b_1(x)y$, and $c=c_2y^2 + c_0(x)$ (with $c_2=c_{02}=\const$).
Also $a\ne 0$ (see the beginning of the proof of Lemma~\lemA).
By (\eqBord) we have
$$
   a\diff{(\Gamma_2)}x + b\diff{(\Gamma_2)}y=ap'+2b_1y^2=(y^2+p)L_1.
$$
Hence $L_1=2b_1$ and then
$$
              ap'=2pb_1.                                       \eqno(\eqGDtwoA)
$$
The second equation in (\eqBord) reads
$$
  b\diff{(\Gamma_2)}x+c\diff{(\Gamma_2)}y=b_1yp' + 2y(c_2y^2+c_0) = (y^2+p)L_2
$$
whence $L_2=2c_2y$ and hence
$c_0 = c_2p - \frac12 b_1 p'$, i.e.,
$c = c_2\Gamma_2 - \frac12 b_1 p'$. By combining this fact with (\eqGDtwoA),
we see that
$$
    c = c_2\Gamma_2 - b_1^2p/a.                                   \eqno(\eqGDtwoB)
$$

If $b=0$, we have $c=c_{02}\Gamma_2$ by (\eqGDtwoB) and $p=\const$ by (\eqGDtwoA).
Then we may set $p=-1$ and we arrive to (iii)--(v) with $n=0$ in (iii).

Now let $b\ne 0$.
The equation (\eqGDtwoA) can be rewritten as
$(\log p)' = 2b_1/a$. Since $p$, $a$, and $b_1$ are polynomials and $\deg(a)\le2$,
we conclude that, up to translation and rescaling the variable $x$,
one of the following three cases occurs.

\medskip
Case 1. $a=1-x^2$.
$$
   \frac{2b_1}{a}=\frac{p'}{p}=\frac{n}{1+x}-\frac{m}{1-x},\qquad
   p=-(1-x)^m(1+x)^n,\qquad m,n>0
$$
(recall that the leading coefficient of $p$ can be chosen arbitrarily).
Therefore $b_1 = \frac12(n-m)-\frac12(n+m)x$. By combining this fact with
(\eqGDtwoB) we obtain (i) as soon as $\Gamma=-\Gamma_2$. So, it remains to show
that $\Gamma$ coincides with $\Gamma_2$ up to a constant factor,
that is $\Gamma_0=\const$.
Indeed, if $\Gamma_0$ were non-constant, by Lemma~\lemChangeW\
it would be a common factor of $a$ and $b$, but
this is impossible for our explicit form of $a$ and $b$.

\medskip
Case 2. $a=x$, $p'/p=2b_1/a=n/x$, $p=x^n$, $n\ge 0$, hence $b_1=n/2$.
Since $a$ and $b$ are coprime, $\Gamma=\Gamma_2$ and we arrive to (ii)
similarly to Case 1. 

\medskip
Case 3. $a=x(x_0-x)$, $p'/p=2b_1/a=-n/(x_0-x)$, $p=-(x_0-x)^n$,
$n\ge 0$, hence $b_1=-nx/2$. This yields (iii).
\qed\enddemo

\proclaim{ Proposition \propDegGaOne }
The following is a complete list of solutions $(a,b,c;\Gamma)$ of
the $(1,\infty)$-AlgDOP problem over $\bC$ up to $(1,\infty)$-admissible change of
variables under condition that $\deg_y\Gamma=1$.
Here $a_{10}$, $a_{20}$, $b_{11}$, $c_{02}$ are constants and
$a_0$, $b_1$, $c_0$, $c_1$ are polynomials in $x$; $\deg a_0\le2$,
$\deg b_1\le 1$.
\roster
\item"(i)"
      $(x^2,\,-nxy,\,n^2y^2-c_0\Gamma_1;x^k\Gamma_1)$,
      $\Gamma_1=x^ny-1$, $n\!\ge\!1$, $k=0,\!1$, $\Delta=-x^2c_0\Gamma_1$;
\item"(ii)"
      $(x^2,\, b_{11}\Gamma-1,\, y^2-c_0\Gamma;\;\Gamma)$, $\Gamma=xy-1$,
      $\Delta=(xy+1-x^2c_0-b_{11}^2\Gamma+2b_{11})\Gamma$;
\item"(iii)"
      $(a_0,\, b_1 y,\, c_{02}y^2+c_1y;\; y)$;
\item"(iv)"
      $(a_{10}x+a_{20}x^2,\, b_{11}xy,\, c_{02}y^2+c_1y;\; xy)$;
\item"(v)"
      $(1-x^2,\, 0,\, c_{02}y^2+c_1y;\; (1-x^2)y)$;
\item"(vi)"
      $(0,1-xy,(1-xy)c_0; 1-xy)$.
\endroster
\endproclaim

We shall see in \S\sectSDOP\ that only (iii)--(v) provide a solution of
the $\ww$-SDOP problem. Moreover, in cases (iv)--(v), this appears
only when $(a,b,c)$ is $(a_{10}x,0,c_{01}y)$ or \hbox{$(1-x^2,0,c_{01}y)$} 
which corresponds to the product of Laguerre polynomials or Jacobi and
Laguerre polynomials.

\demo{ Proof }
Let $\Gamma_1$ be the irreducible factor of $\Gamma$ of $y$-degree $1$, i.e., of the
form $\Gamma_1=\gamma_1 y - \gamma_0$ where $\gamma_0$ and $\gamma_1$ are
polynomials in $x$ and $\gamma_1\ne 0$.

\medskip
Case 1. $a\ne0$ and $\gamma_1$ is not a constant.
By Lemma~\lemValuA, if $x_0$ is a root of $\gamma_1$, it is also a root of $b_1$ and
a double root of $a$.
Since $a$ is a non-zero polynomial in $x$ of degree at most $2$, it has only one
double root, hence $\gamma$ also has only one root (maybe multiple) at $x_0$,
hence we may assume
that $x_0=0$ and
$a=x^2$, $b_1=b_{11}x$, and $\gamma_1=x^n$, $n\ge 1$.
Then $\gamma_0/\gamma_1$ is a
Laurent polynomial, let us denote it by $p(x)=\sum_k p_k x^k$. By a change
of variable $y=y_1+p_0+p_1x+p_2x^2+\dots$ we may kill all the coefficients
of non-negative powers of $x$. Thus we have $p=p_{-n}x^{-n}+\dots+p_{-1}x^{-1}$.
The assumption that $\Gamma_1$ is irreducible implies $\gamma_0(0)\ne 0$, hence
$p_{-n}\ne 0$, i.e., $\ord_x p = -n$.
By rescaling the variable $y$ we may assume $p_{-n}=1$.
The curve $\Gamma_1$ is parametrized by $x=x$, $y=p(x)$. Hence
(\eqParam) yields
$$
     b(x,p)=x^2p',\qquad c(x,p)=b(x,p)p', \qquad\text{ and thus }
      \quad c(x,p)=(xp')^2.                                        \eqno(\eqGDoneA)
$$
Since $b=b_{11}xy+b_0(x)$, the first equation in (\eqGDoneA) reads
$$
   b_{11}\big( x^{-n+1}+\dots +p_{-2}x^{-1}+p_{-1}\big) + b_0 =
             -nx^{-n+1}-\dots-2p_{-2}x^{-1}-p_{-1}.              \eqno(\eqLaurent)
$$

Case 1.1. $n\ge 2$.
Since $b_0$ is a polynomial, we derive from (\eqLaurent) that
$b_{11}=-n$,
$b_0=-p_{-1}-b_{11}p_{-1}=(n-1)p_{-1}$ (a constant),
and $p_{-n+1}=\dots=p_{-2}=0$, i.e., $p=x^{-n}+p_{-1}x^{-1}$.
Then the last equation in (\eqGDoneA) reads
$$
   c_{02}\big(x^{-n}+p_{-1}x^{-1}\big)^2 + (x^{-n}+p_{-1})c_1 + c_0
     = \big(nx^{-n}+p_{-1}x^{-1}\big)^2.                      \eqno(\eqGDoneB)
$$
By comparing the coefficients of $x^{-2n}$ and then those of $x^{-n-1}$,
we obtain $c_{02}=n^2$ and then $p_{-1}=0$. Hence $p=x^{-n}$ and $b_0=0$.
By putting $p_{-1}=0$ into (\eqGDoneB), we get
$n^2 x^{-2n} + x^{-n}c_1 + c_0 = n^2 x^{-2n}$, i.e.,
$c_1=-x^nc_0$. Thus $(a,b,c)=(x^2,\,-nxy,\,n^2y^2-x^nc_0y+c_0)$ which
corresponds to (i) with $k=0$ when $\Gamma=\Gamma_1$.
If $\Gamma$ has another non-constant factor $\Gamma_0$ (a polynomial in $x$),
the condition (\eqBord)
implies that $\Gamma_0$ is a factor of $a$ and $b$
which corresponds to (i) with $k=1$.

\smallskip
Case 1.2. $n=1$ (thus $p=x^{-1}$).
Then (\eqLaurent) reads $b_{11}+b_0=-1$ and 
the last equation in (\eqGDoneA) reads
$c_{02}x^{-2} + x^{-1}c_1 + c_0 = x^{-2}$
whence $b_0=-b_{11}-1$, $c_{02} = 1$, and $c_1=-xc_0$.
Hence $b=b_{11}xy+b_0=b_{11}xy-b_{11}-1$ and
$c = c_{02}y^2+c_1y+c_0=y^2-xyc_0+c_0$ which corresponds to (ii) when $\Gamma=\Gamma_1$.
Once again, any additional non-constant factor $\Gamma_0$ of $\Gamma$ should be a factor
of $a$ and $b$. Then $\Gamma_0=x$ and it divides $b=b_{11}(xy-1)-1$. Therefore
$b_{11}=-1$, thus $b=-xy$. This corresponds to (i) with $n=k=1$.

\medskip
Case 2. $\gamma_1=\const$. Then, up to an admissible change of variables,
we may assume that $\Gamma_1=y$. If $\Gamma=\Gamma_1$, this corresponds to (iii).
Otherwise, as in the previous cases $\Gamma/\Gamma_1$ is
a polynomial in $x$ which divides $a$ and $b$ which corresponds to (iv) and (v).

\medskip
Case 3. $a=0$. Then $\Gamma$ is a factor of $\Delta=-b^2$. By an admissible
change of variables we can reduce $b$ to $y$, $xy$, or $xy-1$.
This leads to (iii), (iv), or (vi).
\qed\enddemo

%%%%%%%%%%%%%%%%%%%%%%%%%%%%%%%%%%%%%%%%%%%%%%%%%%%%%%%%%%%
%%%%%%%%%%%%%%%%%%%%%%%%%%%%%%%%%%%%%%%%%%%%%%%%%%%%%%%%%%%
%%%%%%%%%%%%%%%%%%%%%%%%%%%%%%%%%%%%%%%%%%%%%%%%%%%%%%%%%%%
%%%%%%%%%%%%%%%%%%%%%%%%%%%%%%%%%%%%%%%%%%%%%%%%%%%%%%%%%%%
%%%%%%%%%%%%%%%%%%%%%%%%%%%%%%%%%%%%%%%%%%%%%%%%%%%%%%%%%%%
%%%%%%%%%%%%%%%%%%%%%%%%%%%%%%%%%%%%%%%%%%%%%%%%%%%%%%%%%%%

\head\sectSmallW. Solution of the $(1,w)$-AlgDOP Problem in $\bC^2$ for $1<w\le 2$
\endhead

In this section we find all solutions of the $(1,2)$-AlgDOP Problem in $\bC^2$.
They include all solutions of the $(1,w)$-AlgDOP Problem for any $w$ in the
range $1<w\le 2$ (see \S\sectW).

\subhead \sectHirz. Compactification of $\bC^2$
\endsubhead

Let $(a,b,c;\,\Gamma)$ be a solution of the $(1,2)$-AlgDOP Problem in $\bC^2$,
and let $\Delta=ac-b^2$.
The Newton polygon of $\Delta$ (and hence of $\Gamma$) is contained in the
triangle $[(0,0),(6,0),(0,3)]$ (see Figure~\figABCi).
Therefore it is natural to consider $\bC^2$ as the affine chart $Z\ne 0$
(with coordinates $x=X/Z$, $y=Y/Z^2$)
of the {\it weighted projective plane}
$\bP^2_{1,2,1}$ which is the quotient of $\bC^3\setminus\{(0,0,0)\}$ by the
equivalence relation $(X,Y,Z)\sim(\lambda X,\lambda^2 Y,\lambda Z)$,
$\lambda\ne 0$ (we denote the class of $(X,Y,Z)$ by $[X{:}Y{:}Z]$).
This variety is smooth except at the point $[0\,{:}1{:}\,0]$.

A generic polynomial $P(x,y)$ with $\cN(P)=[(0,0),(6,0),(0,3)]$ defines
an affine curve $\{P=0\}$ in $\bC^2$ whose closure in $\bP^2_{1,2,1}$ is a smooth
curve not passing through the singular point $[0\,{:}1{:}\,0]$, and the linear
projection from this point
$$
    \bP_{1,2,1}^2\setminus\{[0\,{:}1{:}\,0]\}\to\bP^1,\quad
    [X{:}Y{:}Z]\mapsto(X{:}Z),                                  \eqno(\eqWeiPro)
$$
is a 3-fold branched covering of $\{P=0\}$ onto $\bP^1$.
This is why $\bP^2_{1:2:1}$ is an appropriate
compactification of $\bC^2$ in our setting.
However the coefficients of $\Delta$ are not necessarily generic and
the closure of $\{\Gamma=0\}$ may have singularities and
it may pass through $[0\,{:}1{:}\,0]$. To deal with such curves it is convenient
to blow up the point $[0\,{:}1{:}\,0]$. This means that we consider $\bC^2$ as
an affine chart $(x,y)$ of $\cF_2$ -- the {\it Hirzebruch surface of degree} 2 which
is the smooth complex surface obtained by gluing together four copies of $\bC^2$
with coordinates $(x_k,y_k)$, $k=0,\dots,3$, (where $x_0=x$, $y_0=y$ are the coordinates
on $\bC^2$ we started with).
The transition functions are:
$$
  \matrix x_1=1/x,\;\,\qquad & x_2=x,\;\;\;\,\qquad & x_3=\;x_1\;\; = 1/x,\;\;
      \lower3pt\hbox{\mathstrut}  \\
          y_1=y/x^2,  \qquad & y_2=1/y,      \qquad & y_3=1/y_1=x^2/y.
  \endmatrix                                                       \eqno(\eqHirz)
$$
The set of real points of $\cF_2$ is diffeomorphic to a torus.
In Figure~\figHirz\ we represent it as a rectangle with opposite edges identified.
As an illustration of the chart gluing,
we also show in Figure~\figHirz\
how the closures of some two curves in $\cF_2$ look like.

\midinsert
\centerline{\epsfxsize=60mm\epsfbox{hirz.eps}} 
\botcaption{ Figure \figHirz }
   Coordinate axes for all the four charts on $\cF_2$
   and the curves $\{y=1\}=\{y_1=x_1^2\}$ and
   $\{y=-x^4-1\} = \{y_3 = - x_3^2/(1+x_3^4)\}$.
\endcaption
\endinsert

The projection (\eqWeiPro) extends to the projection $\pi:\cF_2\to\bP^1$
given in the affine charts by $(x_k,y_k)\mapsto(x_k:1)$ for $k=0,2$ and
$(x_k,y_k)\mapsto(1:x_k)$ for $k=1,3$. It is a fibration with fiber $\bP^1$.
Its restriction to the closure of a curve $\{P(x,y)=0\}$ is a branched covering
of degree $\deg_y P$.
The strict transform of $[0\,{:}1{:}\,0]$ under the blowup (we denote it by $E$) is
given by $y_2=0$ or $y_3=0$ in the respective charts.
Its self-intersection is $-2$.

The set of $(1,2)$-admissible changes of variables coincides with
the set of biregular automorphisms of $\cF_2$ preserving $\bC^2$.

\subhead\sectIrred. The case when $\Delta$ is irreducible and $\deg_y\Delta=3$
\endsubhead

Let $(a,b,c;\,\Gamma)$ be a solution of the $(1,2)$-AlgDOP Problem such that
$\Gamma=\Delta=ac-b^2$, $\Gamma$ is irreducible, and $\deg_y\Gamma=3$.
We assume that $(a,b,c;\,\Gamma)$ cannot be reduced to a solution of the
$(1,1)$-AlgDOP Problem by a $(1,2)$-admissible change of coordinates.
Let $C$ be the closure of $\{\Gamma=0\}$ in $\cF_2$.
We identify $\bC_2$ (where the affine curve $\Gamma=0$ sits) with the
affine chart corresponding to the coordinate system $(x,y)$.
We shall call it {\it the main chart}.
We denote the fiber $\{x_1=0\}$ by $L_\infty$ (see Figure~\figHirz).
The condition $\deg_y\Gamma=3$ implies that
$C$ is disjoint from $E$, and $\pi|_C$ is a 3-fold branched covering.
Let $\nu:\tilde C\to C$ be the normalization (non-singular model) of $C$.
This means $\tilde C$ is a smooth compact Riemann surface of genus $\bold g$ and $\nu$
a holomorphic mapping which is injective outside a finite number of points.
There is a 1-to-1 correspondence between points of $\tilde C$ and local branches
of $C$.

The genus formula for $C$ reads
$$
   \bold{g}=1+\tfrac12C(C+K_{\cF_2})-\sum_{P\in C}\delta_P
    =4-\sum_{P\in C}\delta_P                                          \eqno(\eqGenusP)
$$
where $\delta_P=\delta_P(C)$ is the delta-invariant of $(C,P)$,
i.e. $2\delta_P=m_i(m_i-1)$ where $m_1,m_2,\dots$ are the multiplicities of
the infinitely near points of $C$ (note that the ``$4$'' in (\eqGenusP) can be
computed as the number of integral points in the interior
of the triangle $[(0,0),(6,0),(0,3)]$).
It is convenient to rewrite (\eqGenusP) in terms of local branches of $C$ as it
is done in [\refBOZ, \S3.2]. Namely,
for a point $P\in C$ with local branches $\gamma_1,\dots,\gamma_r$ we set
$n_P=\sum_{1\le i<j\le r}\gamma_1\cdot\gamma_j$.
Then we have
$\delta_P=n_P+\sum\delta(\gamma_i)$, hence the genus formula (\eqGenusP)
takes the form
$$
  \bold g = 4 - n - \sum_\gamma \delta(\gamma),\qquad
        n = \sum_{P\in C}n_P                                  \eqno(\eqGenus)
$$
where the first sum is over all local branches of $C$.

For a local branch $t\mapsto\gamma(t)=(\xi(t),\eta(t))$
we denote the ramification
index of $\pi\circ\gamma$ by $m_\pi(\gamma)$.
The number $m_\pi(\gamma)$ can be also defined as the intersection
number of $\gamma$ with the fiber of $\pi$ passing through the center of $\gamma$.
If $\ord\xi\ge 0$ then
$m_\pi(\gamma) = \ord_t(\xi(t)-\xi(0))$. If $\ord\xi<0$, then
$m_\pi(\gamma) = -\ord_t\xi$.
By Riemann-Hurwitz formula we have
$$
    2-2\bold g = 6 - \sum_\gamma(m_\pi(\gamma)-1).                 \eqno(\eqRH)
$$

\proclaim{ Lemma \lemPlu}

\noindent(a). Let $\gamma$ be a local branch of $C$ at a point $P$.
Then $m_\pi(\gamma)\le 3$ and we have:
\roster
\item"$\bullet$" 
    if $m_\pi(\gamma)=3$, then $P\in L_\infty$ and $\gamma$ is smooth
\par\noindent
    (we denote the number of such branches by $\beta_3$);
\item"$\bullet$" 
    if $m_\pi(\gamma)=2$ and $\gamma$ is smooth, then $P\in L_\infty$
\par\noindent
    (we denote the number of such branches by $\beta_3$;
         it is clear that $\beta_2+\beta_3\le 1$);
\item"$\bullet$"
    if $\gamma$ is singular, then it is a singularity of type $A_{2k}$
    and $m_\pi(\gamma)=2$
\par\noindent
    (we denote the number of such branches by $\alpha_{2k}$);
\endroster

\smallskip\noindent(b). The curve $C$ is rational (i.e., $\bold g=0$)
and one of the following cases occurs
(among the numbers $n$, $\alpha_k$, $\beta_k$, we list only
the non-zero ones):
\roster
\item"(i)"   $\alpha_2=\alpha_4=\beta_3=n=1$;
\item"(ii)"  $\alpha_2=4$;
\item"(iii)" $\alpha_2=3$ and $\beta_2=n=1$;
\item"(iv)"  $\alpha_2=n=2$ and $\beta_3=1$.
\item"(v)"   $\alpha_4=2$ and $\beta_3=1$;
\item"(vi)"  $\alpha_2=\alpha_6=\beta_3=1$;
\item"(vii)" $\alpha_2=2$ and $\alpha_4=\beta_2=1$.
\endroster

\endproclaim

\demo{ Proof }
(a). Follows from Lemma~\lemValu.
The only point which maybe needs some comments is the smoothness of $\gamma$
in the case when $m_\pi(\gamma)=3$, and hence $P\in L_\infty$.
Up to an admissible change of coordinates we
may assume that $P$ is at $(x_1,y_1)=(0,0)$ (see Figure~\figHirz). Then
$v_\gamma(y_1)>0$ and the condition $m_\pi(\gamma)=3$
means that $v_\gamma(x_1)=3$, i.e.,
$v_\gamma(x)=-3$. Hence Lemma~\lemValu(c) implies that
$v_\gamma(y)\not\in[-4,-1]$. By (\eqHirz) we have
$v_\gamma(y_1)=v_\gamma(y)-2v_\gamma(x)=v_\gamma(y)+6$, thus
$v_\gamma(y_1)\not\in[2,5]$.
It is easy to see that $v_\gamma(y_1)<6$. Hence $v_\gamma(y_1)=1$ whence the result.

\smallskip
(b). We have $\beta_2+\beta_3\le 1$ and the equations (\eqGenus) and (\eqRH) imply
$$
    4 = \bold g + n + \sum_{k\ge 1}k\alpha_{2k}, \qquad
       2\bold g + 4 = \beta_2 + 2\beta_3 + \sum_{k\ge 1}\alpha_{2k}.
$$
The only non-negative solutions are (i)--(vii).
\qed\enddemo

Notice that if a curve in $\bP^2$ is parametrized by
$t\mapsto\big(\xi(t):\eta(t):\zeta(t)\big)$, then the
projectively dual curve is parametrized by
$$
   t\mapsto
          (\dot\eta \zeta - \dot\zeta\eta :
           \dot\zeta\xi   - \dot\xi  \zeta:
           \dot\xi  \eta  - \dot\eta \xi  ).
                \eqno(\eqDual)
$$

\proclaim{ Lemma \lemParam }
The cases (v)--(vii) in Lemma~\lemPlu\ are unrealizable.
In the other cases the curve $C$ admits one of the following parametrizations
$t\mapsto[X:Y:Z]$ in the weighted homogeneous coordinates introduced in
\S\sectHirz:

\roster
\item"(i)" $[\, 32\,(t+1) : 256\,(5t+3)(t+3) : (t+3)^3 \,]$, thus
$$
     \Gamma = y^3 - 20xy^2 + 16y^2 + 45x^3y - 40x^2y - 27x^5 + 25x^4;  \eqno(\eqI)
$$
\item"(ii)" $[\,t^2(t+1) : t^2(2t+1)
                             : 3t+1 + \alpha t^2(t+1)\,]$,
             where  $\alpha^3-9\alpha^2+27\alpha\ne0$;
\vskip3pt
\item"(iii)"  % 3A2 + A1
   $[\,(t-1)^2(t-\alpha) : (t-1)^3(2t^3+t^2-\alpha t^2+t-\alpha t-2\alpha)
                         : (t+\alpha)^2(\alpha t+2t-2\alpha-1)\,]$,
   where $\alpha(\alpha^2-1)(\alpha^2+4\alpha+1)\ne 0$.

\vskip3pt
\item"(iv)" $[\,(t-2)^2(t+1) : (t-2)^3(3t^2 + 3\alpha t + 2\alpha) : 1\,]$,
      where $\alpha\not\in\{-3/2,7/2\}$.
\endroster
\endproclaim

\midinsert
\centerline{\epsfxsize=62mm\epsfbox{blowup.eps}}
\botcaption{Figure \figBlowUp}
\endcaption
\endinsert

\demo{ Proof } In each of the cases (i)--(vii) we may assume that
$C$ is singular at the origin of the main chart.
Then we blow up this point and blow down the strict transform
of the lines $x=0$ and $E$.
In coordinates this means that we consider
the curve $C_1$ on $\bP^2$ which is the projective closure of
the affine curve $\Gamma(x,xy)/x^2=0$. We consider the homogeneous
coordinates $(X_1:Y_1:Z_1)$ on $\bP^2$ such that $x=X_1/Z_1$, $y=Y_1/Z_1$.

Then $C_1$ is a quartic curve tangent to the line $X_1=0$ at $(0{:}1{:}0)$
(see Figure~\figBlowUp).
If $C$ has the $A_{2k}$ singularity at the origin,
then $C_1$ has the $A_{2k-2}$ singularity somewhere on the line $X_1=0$ (when $k>1$)
or a simple tangency (when $k=1$).
The line $Z_1=0$ is the strict transform of $L_\infty$, thus $C_1$
has the tangency with it of the same nature as $C$ has with $L_\infty$.
Up to a $(1,2)$-admissible change of coordinates,
the curve $C$ is determined by $C_1\cup\{X_1Z_1=0\}$.
The weighted homogeneous coordinates (see \S\sectHirz) are expressed
via $(X_1:Y_1:Z_1)$ as follows:
$$
    [X:Y:Z] = [X_1 : X_1 Y_1 : Z_1].                    \eqno(\eqBlowUp)
$$
Now we separately consider the cases of Lemma~\lemPlu.

\medskip
{\it Case} (i). We assume that the node is at the origin.
In this case the curve $C_1$ has singularities $A_2$ and $A_4$.
An irreducible quartic curve with such singularities is unique
up to an automorphism of $\CP^2$ and it is autodual
(see e.g.~[\refBOZ, Cor.~3.10(iii)]).
Such a curve has a single flex point $P$ (because the dual has a single cusp)
and the line $Z_1=0$ is the tangent to $C_1$ at $P$. Let $Q$ be the other
intersection point of $C_1$ with $\{Z_1=0\}$ (see Figure~\figAA). Then
$Q=(0{:}1{:}0)$ and the line $\{X_1=0\}$ is tangent to $C_1$ at $Q$.
This means that $C_1\cup\{X_1Z_1=0\}$ is uniquely determined up to
automorphism of $\bP^2$ whence the uniqueness of $C$.
Thus it remains to check that the curve given
in the statement of the lemma has the required properties.
Indeed, it has a cusp at $[\frac{32}{27}:\frac{256}{81}:1]$ ($t=0$),
an $A_4$ singulariry at $[0{:}0{:}1]$ ($t=\infty$), a node at
$[1{:}1{:}1]$ ($t=-5\pm2\sqrt5$), and a flex point with tangent $L_\infty$ at
$[1{:}0{:}0]$ ($t=-3$).

\midinsert
\centerline{\epsfxsize=50mm\epsfbox{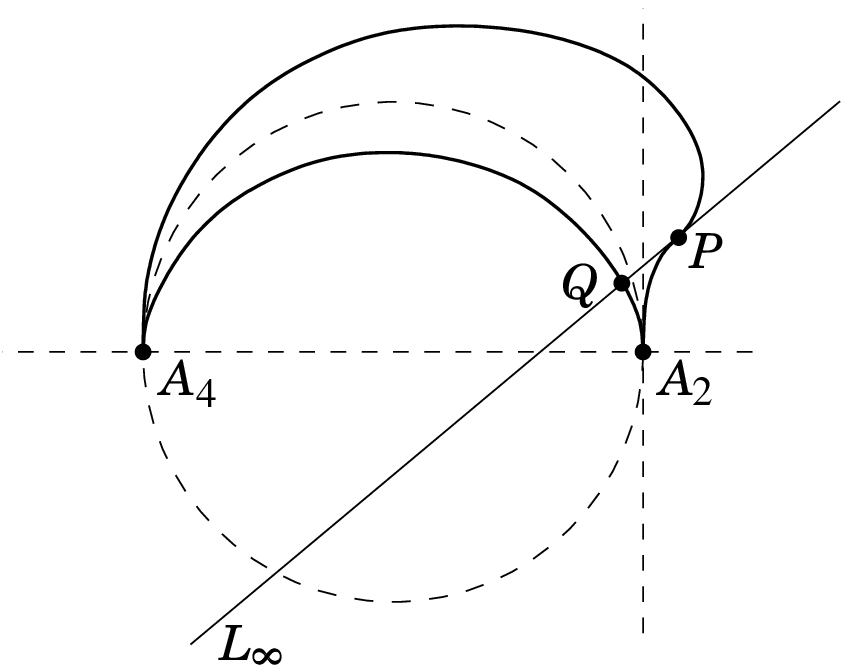}}
\botcaption{ Figure \figAA } A quartic curve with $A_2$ and $A_4$.
\endcaption
\endinsert

\medskip
{\it Case} (ii).
In this case the curve $C_1$ has three cusps $A_2$.
The line $X_1=0$ is its bitangent.
It is well known that such a curve (tricuspidal cubic)
is unique up to automorphism of $\CP^2$, it has only one bitangent line,
and the tangency points $P$ and $Q$
are interchangeable by an automorphism of $(\bP^2,C_1)$.
Thus $C$ is determined by a choice of a line $Z_1=0$ passing through $P$.
So, it depends on one parameter.

The curve $C_1$ is projectively dual to a nodal cubic, hence it has two real forms.
We choose the one with the real bitangent (and hence with two complex conjugate cusps).
In some homogeneous coordinates $(X_2:Y_2:Z_2)$, it has a parametrization
$$
     X_2 = t^2(t+1)^2,\qquad Y_2 = 2t+1,\qquad Z_2 = (t+1)(3t+1).
$$
Here the line $X_2=0$ is bitangent, the cusps correspond to
$t=(-3\pm i\sqrt3)/6$ and $t=\infty$.
According to the above discussion, we may choose $X_1=X_2$, $Y_1=Y_2$, and
$Z_1=Z_2+\alpha X_2$. Plugging this into (\eqBlowUp), we conclude.
The roots of $\alpha^3-9\alpha^2+27\alpha$ are excluded because they correspond
to the cases when the line $L_\infty$ passes through a cusp which
contradicts our assumptions (see Lemma~\lemPlu(a)).

\medskip
{\it Case} (iii).
In this case $C_1$ has two cusps and one node.
We assume that one of the cusps is at the origin.
Then the line $X_1=0$ is bitangent. Let $C_2=\check C_1$
be the curve projectively dual to $C_1$.
Using Pl\"ucker's formulas (in the form given in [\refBOZ, \S3.2] or in [\refKu]),
one can check that the class of curves with two cusps, one node, and at least
one bitangent is stable under the projective duality.
Hence $C_2$ has two cusps and one node. Therefore, in some homogeneous
coordinates $(X_3:Y_3:Z_3)$, it has a parametrization
$$
   X_3 = t^2,\qquad Y_3  =(t-1)(t-\alpha),\qquad Z_3 = t^2(t-1)(t-\alpha).
                                                                    \eqno(\eqIII)
$$
The cusps correspond to $t=0$ and $t=\infty$.
The node corresponds to $t=1$ and $t=\alpha$.
Being dual to $C_2$, the curve $C_1$ is parametrized in some coordinates by
$$
   t\mapsto\varphi(t)=(X_2:Y_2:Z_2)
   =\big( 2(t-1)^2(t-\alpha)^2 : \alpha+1-2t : (\alpha+1)t-2\alpha \big).
$$
(see (\eqDual)).
The points where $C_2$ touches the bitangent line correspond to the local branches
of $C_3$ at the node, which are at $\varphi(1)$ and $\varphi(\alpha)$,
thus $X_1=X_2$. In the coordinates $(X_1{:}Y_1{:}Z_1)$,
the tangency points are at $(0{:}0{:}1)$ and at $(0{:}1{:}0)$.
Up to rescaling the coordinate $t$, we may assume that
$\varphi(1)=(0{:}1{:}0)$. Then the line $Z_1=0$ is uniquely determined by the
condition that it passes through $(0{:}1{:}0)$ and is tangent to $C_1$ which gives
$$
    Z_1 = \tfrac{\alpha}2 X_2 + (\alpha+1)Y_2 + \alpha^2(\alpha+1) Z_2.
$$
The line $Y_1=0$ should pass through $(0{:}0{:}1)$. Then we may set $Y_1=Y_2+Z_2$
and using (\eqBlowUp) we obtain the required parametrization of $C$.

The condition $\alpha\not\in\{0,1\}$ is clear from the construction.
If $\alpha=-1$, then $C_2$ is a double conic because $X_3$, $Y_3$, and $Z_3$
are functions of $t^2$ (see (\eqIII)). If $\alpha$ is a root of
$\alpha^2+4\alpha+1$, then one of the cusps is on $L_\infty$ which
contradicts Lemma~\lemPlu.

\medskip
{\it Case} (iv). The condition $n=2$ can be attained in three different ways.

\midinsert
\centerline{\epsfxsize=120mm\epsfbox{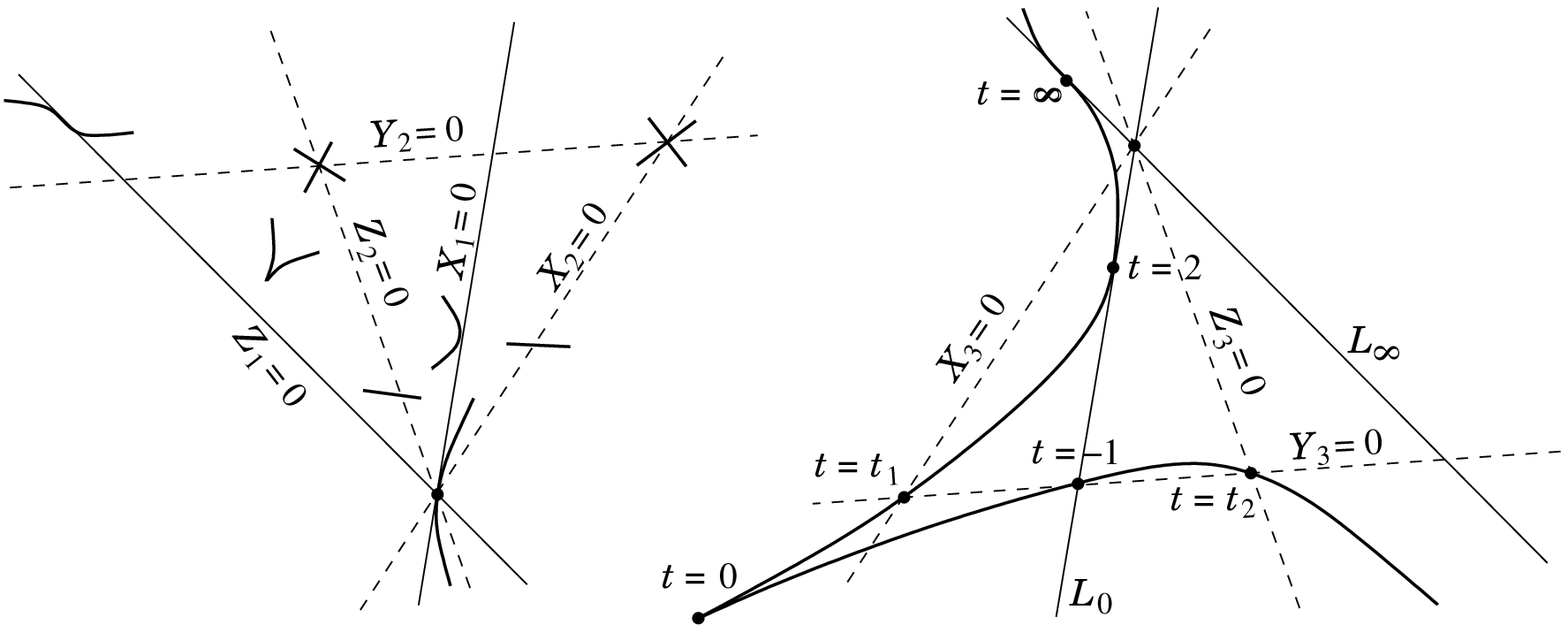}}
\botcaption{ Figure \figIV } Case (iv${}_1$) of Lemma~\lemParam.
\endcaption
\endinsert

\smallskip
{\it Case} (iv${}_1$): $2A_2+2A_1$.
The curve $C$ has two cusps and two nodes.
We assume that one of the cusps is at the origin.
Then $C_1$ has one cusp and two nodes. The line $X_1=0$ is bitangent and the line
$Z_1=0$ is tangent at a flex point. We choose the line $Y_1=0$ to pass through
the both tangency points.
Next we choose
coordinates $(X_2:Y_2:Z_2)$ as shown on the left hand side of Figure~\figIV\
and perform the Cremona transformation
$(X_2:Y_2:Z_2)\mapsto(X_3:Y_3:Z_3)=(Y_2Z_2:Z_2X_2:X_2Y_2)$.
Then the lines
$X_1=0$ and $Z_1=0$ are transformed into lines that we denote by $L_0$ and $L_\infty$.
The transform of $C_1$ is a cuspidal cubic $C_3$ shown
on the right hand side of Figure~\figIV\ (two complex conjugate crossings
with $Z_3=0$ are not shown). The lines $X_2=0$ and
$Z_2=0$ cannot be tangent to the local branches of $C_1$ at the nodes because
otherwise $C_3$ would have too many tangents in the pencil of lines through
$(0:1:0)$ (see Figure~\figIVtan).

\midinsert
\centerline{\epsfxsize=80mm\epsfbox{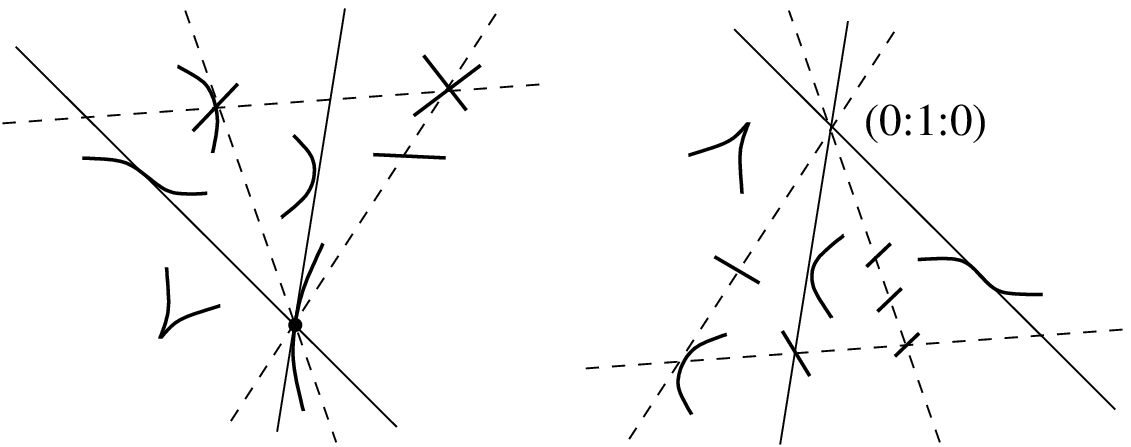}}
\botcaption{ Figure \figIVtan } Unrealizable tangency at $A_1$ in Case (iv${}_1$).
\endcaption
\endinsert

Let $(X_4:Y_4:Z_4)$ be coordinates such
that $C_3$ is parametrized by $t\mapsto\varphi(t)=$ \hbox{$(t^2:t^3:1)$}.
Up to rescaling we may assume
that the point of tangency of $C_3$ with $L_0$ is $\varphi(2)=(4:8:1)$.
Then $C_3\cap L_0=\varphi(-1)=(1:-1:1)$ (see Figure~\figIV) and
the whole configuration is determined by a choice of the line $Y_3=0$
passing through $\varphi(-1)$, i.e., it is determined by a single parameter
$\alpha$ such that $Y_3 = (Y_4+Z_4) - \alpha(X_4 - Z_4)$.
Let $\varphi(t_1)$ and $\varphi(t_2)$ be the other two points of $C_3\cap\{X_3=0\}$.
Then $t_1$ and $t_2$ are the roots of
$(t^3+1) = \alpha(t^2-1)$ different from $-1$, i.e., the roots of
$t^2-(\alpha+1)(t-1)=0$,
hence $t_2=t_1/(t_1-1)$ and $\alpha = (t_1^3+1)/(t_1^2-1)$.
The parametrization of $C_3$ is
$$
   (X_3{:}Y_3{:}Z_3) = \big( t^3 - t_1^3 - 3(t^2 - t_1^2) : (t+1)(t-t_1)(t-t_2) :
                         t^3 - t_2^3 - 3(t^2 - t_2^2) \big)
$$
and after routine computations we obtain the required parametrization of $C$.
If $\alpha\in\{-3/2,7/2\}$, this parametrization defines a curve covered by Case (i).

\smallskip
{\it Case} (iv${}_2$): $A_2+D_5$. The curve $C$ has one cusp and a
singular point of type $D_5$ (given by $u(u^2+v^3)=0$ is some
local curvilinear coordinates). We assume that the $D_5$ singularity
is at the origin.
We may also assume that the line $Y=0$ is tangent to its cuspidal local branch
and passes through the point of tangency of $C$ with the line $Z=0$.
Then $\cN(\Gamma)=[(4,0),(5,0),(0,3),(1,2)]$. Hence after blowing up the origin
we obtain a curve $\bP^2$ with the Newton polygon $[(1,0),(2,0),(0,3),(0,2)]$
(cf.~Figure~\figBlowUp). This is a cuspidal cubic which has a simple (quadratic)
tangency with the line $X_1=0$ and a cubic tangency with the line $Z_1=0$.
One easily checks that these conditions uniquely determine $C\cup\{X_1Z_1=0\}$
up to automorphism of $\bP^2$. Hence the curve $C$ is unique up to
an admissible change of coordinates. It remains to observe that the parametrization
in the statement of the lemma gives the required curve when $\alpha=-1$.
It has a cusp at $[0:0:1]$ and the $D_5$ singularity at $[4:16:1]$.

\smallskip
{\it Case} (iv${}_3$): $2A_2+A_3$. The curve $C$ has two cusps $A_2$ and one
tacnode $A_3$ (ordinary tangency of two smooth branches).
We assume that one of the cusps is at the origin.
Then the curve $C_1$ has one cusp $A_2$, one tacnode $A_3$ and at least one
bitangent (the line $X_1=0$). Let us show that this combination is impossible
for quartic curves in $\bP^2$. It is convenient to apply Pl\"ucker's equations
in the form given in [\refKu, Thm.~1.3]. In the notation of [\refKu] we
have $d=4$, $g=0$, $n_v=2$, $c_v=1$, hence $\hat d=5$ by [\refKu, Eq.~(1.6)].
Then the equations [\refKu, (1.8)--(1.9)] read
$\hat n+\hat c=6$ and $2\hat n+3\hat c=16$ whence $\hat n=2$.

It is easy to check that the dual branches of a tacnode form a tacnode of
the dual curve and it contributes $2$ to $\hat n$. Hence the dual curve does not
have nodes, i.e., $C_1$ does not have bitangents.
A contradiction.

\medskip
{\it Case} (v). In this case $C_1$ has singularities $A_2$ (on the line $X_1=0$)
and $A_4$. An irreducible quartic curve with such singularities is unique
up to an automorphism of $\CP^2$ and it is autodual (cf.~Case (i)).
Hence it has a unique flex point because the dual curve has one cusp $A_2$.
By the uniqueness, we may assume that $C_1$ is a small real
perturbation of real double conic as explained in the remark in the
proof of [\refBOZ, Cor.~3.10] (see Figure~\figAA).
The line $Z_1=0$ is the tangent at the flex point $P$. The line $X_1=0$ is
the tangent at the transverse crossing $Q$ of $Z_1=0$ with $C_1$ and it must pass
through the cusp $A_2$.
We see in Figure~\figAA\ that this is impossible because the arc $A_2QA_4$
is convex.

\medskip
{\it Case} (vi). We assume that the singularity of $C$ at the origin is $A_6$.
Then $C_1$ has $A_2$ and $A_4$ as singularities and the proof is the same as in
Case (v) but with $A_2$ and $A_4$ exchanged.

\medskip
{\it Case} (vii). We assume that the singularity of $C$
at the origin is $A_2$. Then $X_1=0$ is a bitangent and
$C_1$ has singular points $A_2$ and $A_4$. As we mentioned in Case (v), such
a curve is autodual, hence it cannot have
a bitangent because the dual curve does not have nodes.
A contradiction.
\qed\enddemo

\proclaim{ Proposition \propIrred } Let $(\g;\,\Gamma)$ be
a solution of the $(1,2)$-AlgDOP problem over $\bR$
such that $\Gamma$ is irreducible and $\deg_y\Gamma=3$.
Then, up to a $(1,2)$-admissible change of variables,
either $(\g;\,\Gamma)$ is a solution of the $(1,1)$-AlgDOP problem, or
$\Gamma$ is given by (\eqI) and
$$
g=\left(\matrix
  y + 8x - 9x^2     &  5(4y - 3xy - x^2)        \lower2pt\hbox{\mathstrut}\\
  5(4y - 3xy - x^2) &  -25(y^2 - 4xy + 3x^3) \endmatrix\right).
   \eqno(\eqIG)
$$
\endproclaim

\demo{ Proof }
Given a parametrization $(X(t),Y(t),Z(t))$ in the weighted homogeneous
coordinates introduced in \S\sectHirz, the relation (\eqParam)
applied to $\xi(t)=X(t)/Z(t)$, $\eta(t)=Y(t)/Z(t)^2$ 
yields a system of homogeneous linear equations 
for the coefficients $a_{ij}$, $b_{ij}$, $c_{ij}$ of the entries of the cometric $g$
in each of the cases (i)--(iv) of Lemma~\lemParam.
Up to a constant factor, the only non-zero solution in Case (i) is (\eqIG).

In Cases (ii)--(iv) the system of equations polynomially depends on $\alpha$.
It has a non-zero solution if and only if the gcd of the determinants of
the maximal minors vanishes. A straightforward computation shows
that this gcd is a polynomial in $\alpha$ all whose roots are excluded
in Lemma~\lemParam. For example, in Case (ii), the gcd is a power of $\alpha$
multiplied by a power of $\alpha^2-9\alpha+27$.
\qed\enddemo

\medskip\noindent
{\bf Remark \remIrred.}
The solution of the $(1,2)$-DOP Problem given in [\refBB, \S8]
transforms (up to a constant factor) into our solution given
in Proposition~\propIrred\ by the change of variables
$$
   \theta_1 = \frac{2+\sqrt5}{5}(x-1), \qquad
   \theta_2 = -\frac{\sqrt5}{125}(4y-5x+1).
$$

%%%%%%%%%%%%%%%%%%%%%%%%%%%%%%%%%%%%%%%%%%%%%%%%%%%%%%%%%%%%%%%%%
%%%%%%%%%%%%%%%%%%%%%%%%%%%%%%%%%%%%%%%%%%%%%%%%%%%%%%%%%%%%%%%%%
%%%%%%%%%%%%%%%%%%%%%%%%%%%%%%%%%%%%%%%%%%%%%%%%%%%%%%%%%%%%%%%%%
%%%%%%%%%%%%%%%%%%%%%%%%%%%%%%%%%%%%%%%%%%%%%%%%%%%%%%%%%%%%%%%%%

\subhead \sectDegTwo. The case when $\Gamma$ has a factor of $y$-degree $2$
\endsubhead

\proclaim{ Proposition \propReducTwo }
Let $(\g;\,\Gamma)$ be a solution of the $(1,2)$-AlgDOP Problem
such that $\deg_y\Gamma=2$. Suppose that it is not a solution of
the $(1,\infty)$-AlgDOP problem.
Then, by a $(1,2)$-admissible change of variables,
$(\g;\,\Gamma)$ can be reduced either to a solution
of the $(1,1)$-AlgDOP problem or to one of the following three cases:

\roster
\item"(i)" $\Gamma = y^2-x^3$ and
$$
   g = \left(\matrix 4y & 6x^2 \\ 6x^2 & 9xy+\alpha\Gamma\endmatrix\right)
       + (\beta x+\mu)\left(\matrix 4x & 6y \\ 6y & 9x^2\endmatrix\right).
                                                                       \eqno(\eqReduc)
$$
Rescaling the coordinates, one can replace
 $(\alpha,\beta,\mu)$ by
 $(\lambda\alpha,\lambda\beta,\lambda^{-1}\mu)$ for any non-zero $\lambda$.

\medskip
\item"(ii)" $\Gamma=y(y-x^2)$, $\g=\g_{(\alpha,\beta,\mu)}$
with $(\alpha,\beta-\beta^2,\mu)\ne(0,0,0)$, $\mu\in\{0,1\}$, where
$$
   g_{(\alpha,\beta,\mu)} = (y-x^2)\left(\matrix 1&0\\0&\alpha y\endmatrix\right)
      +(\beta x+\mu)\left(\matrix x&2y\\2y& 4xy\endmatrix\right).
$$

\medskip
\item"(iii)" $\Gamma=y(y-x^2+1)$ and $\g=\g_{(\alpha,\beta)}$
with $(\alpha,\beta-\beta^2)\ne(0,0)$,
where
$$
   \g_{(\alpha,\beta)}
       = (y-x^2+1)\left(\matrix 1&0 \\ 0&\alpha y\endmatrix\right)
       + \beta\left(\matrix x^2-1 & 2xy\\2xy & 4x^2y\endmatrix\right).
                                                            \eqno(\eqReducIII)
$$
\endroster
\endproclaim

\demo{ Proof }
Lemmas~\lemChangeW\ and \lemValuA\ imply that
$\Gamma$ does not have monomials of the form $x^ky^2$ with $k>0$.
Local branches of $\Gamma$ correspond to edges of $\cN(\Gamma)$ (see Lemma~\lemNewton).
Hence Lemma~\lemValu(d) implies that the slope of any upper edge is steeper
than 1:2, hence $\cN(\Gamma)$ is contained in the triangle $[(0,0),(4,0),(0,2)]$.
This fact combined with Lemma~\lemValuA\
(which means that the affine curve $\Gamma=0$ has no vertical tangent)
leaves only five possibilities for $\Gamma$ up to admissible change of coordinates:
the three cases (i)--(iii) and also $y(y-x)$ and $y^2-1$.

In each case, the condition (\eqParam) yields a system homogeneous linear
equations for the coefficients of $a$, $b$, and $c$ (the entries of $g$).
By solving these systems we obtain
the result. In the last two cases we obtain $a_1=0$, thus these are solutions
of the $(1,\infty)$-AlgDOP problem (see Figure~\figABCii).
In the other cases we normalize the solutions so that $a_1=1$.
The parameter $\mu$ in Case (ii) can be set to $0$ or $1$
by rescaling the coordinates (see Example~\exaChangeVarII). The conditions
$(\alpha,\beta-\beta^2,\mu)\ne(0,0,0)$ (in Case (ii)) and
$(\alpha,\beta-\beta^2)\ne(0,0)$ (in Case (iii)) are equivalent to $\det g\ne0$.
\qed\enddemo

\proclaim{ Lemma \lemReduc } Let $(\g;\,\Gamma)$ be a solution of
$(1,2)$-AlgDOP problem which cannot be reduced to a solution of
$(1,1)$-AlgDOP problem.
Let $\Gamma=y(y-p_1(x))(y-p_2(x))$. Then the polynomial $p_1p_2$
has at most two roots (maybe multiple).
\endproclaim

\demo{ Proof } Proposition~\propReducTwo\ applied to $(\g,y(y-p_k))$
implies that $a_0$ vanishes at the roots of $p_k$ (recall that
$g^{11}=a_{01}y+a_0(x)$). It remains to prove that $a_0$ cannot be
identically zero. Indeed, if it is, then $y$ would divide each entry of $\g$,
hence $y^2$ would divide $\det\g$ which is impossible because $\det\g=\Gamma$
in our case (see Figure~\figABCi) and $\Gamma$ is squarefree.
\enddemo

\proclaim{ Proposition \propReducThree }
Let $(\g;\,\Gamma)$ be a solution of the $(1,2)$-AlgDOP problem over $\bC$
such that $\Gamma$ is reducible and $\deg_y\Gamma=3$.
Then, up to a $(1,2)$-admissible change of variables,
either $(\g;\,\Gamma)$ is a solution of the $(1,1)$-AlgDOP problem, or
$g$ is given by (\eqReduc) with $(\alpha,\beta,\mu)=(-18,-3/2,1/2)$ and
hence $\Gamma=\Gamma_1\Gamma_2$ where $\Gamma_2=y^2-x^3$ and
$\Gamma_1=8y - 3x^2 - 6x + 1$.

In the latter case the curve $\Gamma=0$ has singularities of the types
$A_1$, $A_2$, and $A_5$ at $(\frac{1}{9},-\frac{1}{27})$, $(0,0)$, and
$(1,1)$ respectively.
\endproclaim

\demo{ Proof }
We have $\Gamma=\Gamma_1\Gamma_2$ avec $\deg_y\Gamma_k=k$.
By Proposition~\propGammaOne, $(\g,\Gamma_2)$ is also a solution of the $(1,2)$-AlgDOP
problem. Moreover, $(\g,\Gamma)$ reduces to a solution of a $(1,w)$-AlgDOP
problem by a $(1,2)$-admissible change if and only if the same
is true for $(\g,\Gamma_2)$.
Hence we may assume that $(\g,\Gamma_2)$, 
is as in Proposition~\propReducTwo.

\medskip
{\it Case 1.} $\Gamma_2$ is irreducible. Then $(\g,\Gamma_2)$
is as in Proposition~\propReducTwo(i). A computation shows that
$\det\g=\Gamma_1\Gamma_2$ with $\Gamma_1=\alpha y-f(x)$, hence $\alpha\ne 0$.
Then, for the parametrization $\xi=t$, $\eta=f(t)/\alpha$ of $\{\Gamma_1=0\}$,
the equations (\eqParam) take the form
$$
  a(\xi,\eta)\dot\eta-b(\xi,\eta)\dot\xi=6FG/\alpha^2=0, \qquad
  b(\xi,\eta)\dot\eta-c(\xi,\eta)\dot\xi=FGH/\alpha^2=0,
$$
where $F=At-3B$,
$A=(\alpha - 12\beta)(\alpha - 9\beta)$, $B=18 + 5\alpha\mu - 36\beta\mu$,
$G=(\beta t+\mu)^2-t$, and
$H=(9\beta - \alpha)t + 9\mu$.
Since $G$ cannot vanish identically, we have $A=B=0$.

If $\alpha=9\beta$, then $B=-27(2 + \beta\mu)$, and we obtain
a solution of the $(1,1)$-AlgDOP problem (the one discussed in [\refBOZ, \S4.10]).

If $\alpha=12\beta$, then $B=-18 (3 + 4\beta\mu)$, and we obtain the
announced solution.

\medskip
{\it Case 2.} $\Gamma_2$ is reducible. In this case
$\{\Gamma=0\}=L_1\cup L_2\cup L_3$ where $L_k=\{y=p_k(x)\}$, $k=1,2,3$.
Let
$P_1=L_1\cap(L_2\cup L_3)$,
$P_2=L_2\cap(L_3\cup L_1)$, and
$P_3=L_3\cap(L_1\cup L_2)$.
Then Lemma~\lemReduc\ implies that each $P_k$ has at most two points.
By Proposition~\propReducTwo, if $k\ne m$, then $L_k$ and $L_m$
either are tangent, or cross at two points.
Hence, up to a $(1,2)$-admissible change of coordinates,
$\Gamma=y(y-p(x))(y-\lambda p(x))$ where $\lambda\not\in\{0,1\}$
and $p(x)$ is $x^2$ or $x^2-1$.
In this case $\g$ is as in Proposition~\propReducTwo\ (ii) or (iii).
Then one easily checks that (\eqParam) is not satisfied for the parametrization
$\xi(t)=t$, $\eta(t)=\lambda p(t)$ of the curve $y=\lambda p(x)$.
\qed\enddemo

%%%%%%%%%%%%%%%%%%%%%%%%%%%%%%%%%%%%%%%%%%%%%%%%%%%%%%%%%%%%%%%%%
%%%%%%%%%%%%%%%%%%%%%%%%%%%%%%%%%%%%%%%%%%%%%%%%%%%%%%%%%%%%%%%%%
%%%%%%%%%%%%%%%%%%%%%%%%%%%%%%%%%%%%%%%%%%%%%%%%%%%%%%%%%%%%%%%%%
%%%%%%%%%%%%%%%%%%%%%%%%%%%%%%%%%%%%%%%%%%%%%%%%%%%%%%%%%%%%%%%%%
%%%%%%%%%%%%%%%%%%%%%%%%%%%%%%%%%%%%%%%%%%%%%%%%%%%%%%%%%%%%%%%%%
%%%%%%%%%%%%%%%%%%%%%%%%%%%%%%%%%%%%%%%%%%%%%%%%%%%%%%%%%%%%%%%%%

\head\sectSDOP. Solution of the weighted DOP/SDOP problem in $\bR^2$
\endhead

\subhead\sectCompact. Compact solutions
\endsubhead

In what follows we give the measure density $\rho$ without the normalizing
constant which is always assumed to be equal to $1/\int_\Omega\rho\,dx$.

\proclaim{ Theorem \thCompTwo }
Let $(\Omega,g,\rho)$ be a solution of the $(1,2)$-DOP problem in $\bR^2$
which is not a solution of the $(1,\infty)$-DOP problem
and such that $\Omega$ is bounded.
Then, up to a $(1,2)$-admissible change of variables, either it is
a solution of the $(1,1)$-DOP problem, or one of the following cases occurs
(see Figure~\figCompTwo):
\roster
%===================================
  \item"(\modelBi)"
         (dodecahedral quotient)
         $\g$ is given by (\eqIG), $\Gamma=-\frac1{25}\det\g$ is given by (\eqI),
         $\Omega$ is the bounded component of $\bR^2\setminus\{\Gamma=0\}$,
         and $\rho=\Gamma^{p-1}$ with $p>\frac3{10}$;

%===================================
\smallskip
  \item"(\modelBii)"
      (cuspidal cubic with cubically tangent parabola)
      $g$ is as in Prop.~\propReducThree, i.e.,
$$
    \g=g_{(-18,-\frac32,\frac12)}
       =\left(\matrix 4(2y-3x^2+x)  && 6(y-3xy+2x^2)\\
                      6(y-3xy+2x^2) && 9(x^2+x^3+2xy-4y^2)\endmatrix\right),
$$
    $\frac1{36}\det\g=\Gamma=\Gamma_1\Gamma_2$ with
    $\Gamma_1=8y - 3x^2 - 6x + 1$, $\Gamma_2=x^3-y^2$, the domain
    $\Omega$ is the bounded component of $\bR^2\setminus\{\Gamma=0\}$,
    and $\rho=\Gamma_1^{p-1}\Gamma_2^{q-1}$ with $p>0$, $q>\frac16$, $p+q>\frac23$;

%===================================
\smallskip
  \item"(\modelBiii)"
    (parabolic biangle)
    $\g=\g_{(\alpha,\beta)}$ is given by (\eqReducIII) 
    with $\alpha<0$ and $\beta\le 0$;
    $\Omega=\{x^2-1<y<0\}$,
    and $\rho=(-y)^{p-1}(y-x^2+1)^{q-1}$ with $p,q>0$.
\endroster

\smallskip\noindent
The solution (\modelBiii) reduces to a solution of the
$(1,1)$-problem by a $(1,2)$-admissible change
if and only if either $\alpha=4\beta$ (then it is already so),
or $\alpha=4\beta-4$. In the latter case, the variable change is
$(x,y)\mapsto(x,x^2-y-1)$ which transforms $\g_{(4\beta-4,\beta)}$ into
$-\g_{(4-4\beta,1-\beta)}$.
We have $\g_{(4\beta,\beta)}=-\beta G'_{-1/\beta}$ 
in the notation of [\refBOZ, \S4.5].
\endproclaim

\midinsert
\centerline{
  \hbox{\noindent\hskip0pt
   \epsfysize=5cm\epsfbox{a1a2a4.eps}\hskip5mm
   \epsfysize=5cm\epsfbox{a1a2a5.eps}
   \hskip-109.5mm
   \lower-3mm\hbox{$(0,0)$}\hskip41mm
   \lower-13mm\hbox{$(1,1)$}\hskip2pt
   \lower-46mm\hbox{$\big(\frac{32}{27},\frac{256}{81}\big)$}\hskip-8mm
   \lower-11mm\hbox{$(0,0)$}\hskip6mm
   \lower-7mm\hbox{$\big(\frac19,-\frac1{27}\big)$}\hskip8.5mm
   \lower-47mm\hbox{$(1,1)$}\hskip3mm
  }
}
\centerline{Dodecahedral quotient \hskip19mm Cubic $y^2=x^3$ and a parabola}
\botcaption{Figure \figCompTwo}
      The first two domains in
      Theorem~\thCompTwo.
\endcaption
\endinsert

\demo{ Proof }
According to Propositions~\propIrred, \propReducTwo, and \propReducThree,
these are the only solutions of the $(1,\infty)$-AlgDOP problem with
bounded components of $\bR^2\setminus\{\Gamma=0\}$.

The imposed restrictions on $\alpha$ and $\beta$ in Case (\modelBiii)
are equivalent to the positive definiteness of $\g$.
Indeed, we have $\g(0,y)=\diag\big(y-\beta+1,\alpha y(1+y)\big)$, thus
the positivity of $\g$ on $\Omega\cap\{x=0\}$ implies $\alpha<0$ and $\beta\le 0$.
Conversely, let $\alpha<0$ and $\beta\le 0$.
Then $\g>0$ at $(0,-\frac12)$, hence
it is enough to show that $\Delta$ does not vanish in $\Omega$.
We have $\Gamma=y\Gamma_1$ and $\Delta:=\det(\g)=y\Gamma_0\Gamma_1$ where
$$
   \Gamma_0=\alpha y + (4\beta-\alpha)(1-\beta)x^2+\alpha(1-\beta),
   \qquad \Gamma_1=y-x^2+1.                                         \eqno(\eqCompInf)
$$
Thus it is enough to show that $\{\Gamma_0=0\}\cap\Omega=\varnothing$.
If $\beta=0$, then $\Gamma_0=\alpha\Gamma_1$ and we are done.
If $\beta>0$, it is easy to check that the curve $\Gamma_0=0$ does not
cross $\partial\Omega$.

The form of the measure density follows from Proposition~\propMes\ unless
$\Delta$ has a multiple factor. This happens only in Case (\modelBiii) and,
as one can see from (\eqCompInf), only when $\beta=0$
(then $\Delta=\alpha y\Gamma_1^2$) or $\beta=1$ (then $\Delta=\alpha y^2\Gamma_1$).
In these two cases one can perform the computations described in the beginning
of \S\sectALGtoDOP\ (in fact, one case is reduced to the other one by the variable
change indicated in the statement of this theorem). The inequalities for
$p$ and $q$ are the integrability conditions (see [\refBOZ, Remark~2.28]).
\qed\enddemo

\proclaim{ Theorem \thCompInf }
Let $(\Omega,g,\rho)$ be a solution of the $(1,\infty)$-DOP problem in $\bR^2$
such that $\Omega$ is bounded. Then one of the following cases occurs
up to a $(1,\infty)$-admissible change of variables.

\medskip\noindent
\roster
\item"(\modelBiv)"
         $(g,\Gamma)$ is as in Proposition~\propDegGaTwo(i) with $c_{02}<0$,
         $\Omega=\{\Gamma>0\}\cap\{x^2<1\}$ which is the only bounded
         component of $\bR^2\setminus\{\Gamma=0\}$, and $\rho$ is as follows:

\smallskip\noindent
         \hbox to 100mm{\hskip-10pt\vbox{
         \roster
         \item"$\bullet$"
           if one of $m,n$ is odd, then
           $\rho=\Gamma^{p-1}$ with
           $p>\max\big(\frac12-\frac1n,\frac12-\frac1m\big)$;
         \smallskip
         \item"$\bullet$"
           if $m$ and $n$ are both even, then
           $$
             \rho=
             \big( (1-x)^{\frac{m}2}(1+x)^{\frac{n}2} + y\big)^{p-1}
             \big( (1-x)^{\frac{m}2}(1+x)^{\frac{n}2} - y\big)^{q-1}
           $$
           with positive $p$ and $q$ such that $p+q>\max\big(1-\frac2n,1-\frac2m\big)$.
         \endroster
         }} % end of \hbox{\vbox

\smallskip
\item"(\modelBv)"
         $(g,\Gamma)$ is as in Proposition~\propDegGaTwo(iii)
         with $k=x_0=1$, $n\ge 1$, and $c_{02}\le 0$,
         (i.e., $\Gamma=x\Gamma_2$, $\Gamma_2=(1-x)^n-y^2$),
         $\Omega$ is $\{\Gamma>0\}\cap\{0<x<1\}$ which is the only bounded
         component of $\bR^2\setminus\{\Gamma=0\}$, and $\rho$ is as follows:

\smallskip\noindent
         \hbox{\hskip-10pt\vbox{
         \roster
         \item"$\bullet$"
           if $n$ is odd, then
           $\rho=x^{r-1}\Gamma_2^{p-1}$ with $r>0$ and
           $p>\max\big(0,\frac12-\frac1n\big)$;
         \smallskip
         \item"$\bullet$"
           if $n$ even, then
           $\rho=x^{r-1}
             \big( (1-x)^{n/2} + y\big)^{p-1}
             \big( (1-x)^{n/2} - y\big)^{q-1}$
           with \par\noindent
           positive $p$, $q$, and $r$ such that $p+q>1-\frac2n$.
         \endroster
         }} % end of \hbox{\vbox

\endroster
\endproclaim

\demo{ Proof }
According to Proposition~\propDegGaTwo,
the only solutions of the $(1,\infty)$-AlgDOP problem with
bounded components of $\bR^2\setminus\{\Gamma=0\}$ are the indicated ones.

\medskip
Case (\modelBiv). We have
$$
   \Delta:=\det(g) = \Gamma_0\Gamma, \qquad
   \Gamma_0 = \tfrac14\Big((n-m)-(n+m)x\Big)^2 - c_{02}(1-x^2).
$$
The condition $c_{02}<0$ is equivalent to the fact that $\g$ is positive definite.
Indeed, suppose that $c_{02}\ge0$. Then
$$
   \Gamma_0(x_0)=-\frac{4c_{02}(c_{02}+mn)}{4c_{02}+(m+n)^2}\le 0
   \qquad\text{for}\quad
   x_0 = \frac{n^2-m^2}{4c_{02}+(m+n)^2}.
$$
Since $|x_0|<1$, we have $(x_0,0)\in\Omega$ and hence $\Gamma(x_0,0)>0$.
Therefore $\Delta(x_0,0)=\Gamma_0(x_0)\Gamma(x_0,0)\le 0$, thus $\g$
is not positive definite on $\Omega$. Conversely, if $c_{02}<0$, then
$\Gamma_0(x)\ge-c_{02}(1-x^2)>0$ when $|x|<1$, whence $\Delta|_\Omega>0$
which implies $\g|_\Omega>0$ by Sylvester's criterion because $a|_\Omega>0$.

The required form of $\rho$ can be derived from Proposition~\propMes.
Indeed, $\Gamma$ is maximal for $\g$ (by Proposition~\propDegGaTwo) and
$\Delta$ is squarefree. Hence, by Proposition~\propMes,
$\rho$ is of the required form but with an additional factor $\exp Q$.
Let us show that $Q$ is constant.
Set $l=\deg c_0(x)$. Then $(\g,\Gamma)$ is a solution of $(2,l)$-SDOP problem.
We have $c_0(x) = (1-x)^{m-1}(1+x)^{n-1}\Gamma_0$, hence
$l=m+n-2+\deg_x\Gamma_0$ and we obtain $l=\deg_x\Delta-2$.
Thus (\eqMesIneq) with $\ww=(2,l)$ implies $\deg_x Q=0$.
On the other hand, $(\g,\Gamma)$ is a solution of $(1,w)$-SDOP problem
for some $w>2$, hence (\eqMesIneq) for these weights implies $\deg_y Q=0$.
Thus $Q=\const$, hence $\rho$ is of the required form. The inequalities for
$p$ and $q$ are the integrability conditions (see [\refBOZ, Remark~2.28]).

\medskip
Case (\modelBv).
The proof is almost the same as in Case (\modelBiv).
We have $\Delta:=\deg\g=\Gamma_0\Gamma$ where $\Gamma_0=\tfrac14 n^2x-c_{02}(1-x)$.
The condition $c_{02}\le 0$ is equivalent to $\g>0$.
Indeed, if $c_{02}>0$, then $\Gamma(x_0)=0$ and $0<x_0<1$
for $x_0=c_{02}/(c_{02}+\frac14n^2)$, which implies that $\Delta(x_0,0)=0$
whence $\g|_\Omega$ is not positive definite.
Conversely, if $c_{02}\le 0$, then
$\Gamma_0(x)>0$ when $0<x<1$, hence $\Delta|_\Omega>0$ whence
$\g|_\Omega>0$ (because $a|_\Omega>0$).

The form of $\rho$ is established by the same arguments as in Case (\modelBiv). 
\qed\enddemo

%%%%%%%%%%%%%%%%%%%%%%%%%%%%%%%%%%%%%%%%%%%%%%%%%%%%%%%%%%%%%%%%%
%%%%%%%%%%%%%%%%%%%%%%%%%%%%%%%%%%%%%%%%%%%%%%%%%%%%%%%%%%%%%%%%%

\subhead\sectNonComp. Non-compact solutions
\endsubhead

\proclaim{ Theorem \thNonCTwo }
Any solution $(\Omega,g,\rho)$ of the $(1,2)$-SDOP problem with an
unbounded domain $\Omega$ can be
reduced to a solution the $(1,1)$- or $(1,\infty)$-SDOP problem by
a $(1,2)$-admissible change of variables.
\endproclaim

\demo{ Proof } Let $\Delta=\deg\g$ and $\Gamma$ be the maximal boundary for $\g$
(see Definition~\defMaxBndry).
The integrability condition for the measure $\rho\,dx$ implies
that the term $\exp(Q)$ in Proposition~\propMes\ is non-constant,
hence $\deg_y\Gamma\le\deg_y\Delta\le 2$ (see Corollary~\corMes).

\smallskip
Case 1. $\deg_y\Gamma=2$. Then $\g$ is as in Proposition~\propReducTwo.
The condition $\deg_y\Delta<3$ implies $\alpha=0$ in all the three cases (i)--(iii).
Then, using the algorithm of \S\sectALGtoDOP, we find that there is no
exponential term in $\rho$.

\smallskip
Case 2. $\deg_y\Gamma=1$. Then $\deg_y a=1$ and
$\Gamma=y$ up to a $(1,2)$-admissible coordinate change,
because otherwise we have a solution of the $(1,\infty)$-SDOP problem
(see Lemmas~\lemChangeW--\lemValuA\ and Figure~\figABCii).
The fact that $\Gamma=y$ combined with (\eqParam) implies that $y$ divides $b$ and $c$.
Since $\deg_y a=1$ and $\deg_y b\le 1$,
the condition $\deg_y(ac-b^2)<3$ implies $\deg_y c=1$,
hence $(a,b,c)=(y+a_0,yb_1,yc_1)$ and $\Delta = (c_1-b_1^2)y^2 + a_0c_1 y$.
If $\deg c_1\le 1$, then we have a solution of the $(1,1)$-SDOP problem,
hence $\deg c_1=2$. 

\medskip
Case 2.1. $a_0\ne0$. Then $y^2$ does not divide $\Delta$, hence,
by Proposition~\propMes, $\rho=e^h y^\alpha$ for a polynomial $h$ such that
$\deg_y h\le 1$, thus $h=yh_1(x)+h_0(x)$.
By (\eqDh) we have $L^k\in\cP_\ww(w_k)$, $\ww=(1,2)$, where
$$
     L^1=(y+a_0)(yh'_1+h'_0)+yb_1h_1,  \qquad
     L^2=yb_1(yh'_1+h'_0)+yc_1h_1.
$$
The condition $L^1\in\cP_\ww(w_1)$ means that $\deg_y L^1=0$, $\deg_x L^1\le 1$.
Hence $h'_1=0$ (thus $h_1$ is a constant) and $h'_0=-h_1b_1$.
The integrability condition implies that $h_1<0$
(we assume here that $\Omega=\{y>0\}$). Up to rescaling $y$ we may assume that $h_1=-1$,
hence $h_0'=b_1$.
Then we have $L^1=a_0b_1$ and $L^2=(b_1^2-c_1)y$. The condition $L^2\in\cP_\ww(w_2)$
implies that $b_1^2-c_1$ is a constant
and hence $\deg b_1=1$ (recall that we have assumed that $\deg c_1=2$).
Since $L^1=a_0b_1$ and $\deg_x L^1\le 1$, we conclude that $a_0$ is constant.
Translating $x$, we may achieve $b_{01}=0$ and we obtain
$(a,b,c)=(y+\alpha,\beta xy,\beta^2x^2y+\gamma y)$ for some constants
$\alpha,\beta,\gamma$. The change $(x,y)\mapsto(x,y-\frac12\beta x^2)$
transforms it into
$(y+\frac12\beta x^2+\alpha,-\alpha\beta x,\alpha\beta^2x^2+\gamma y)$
which is a solution of the $(1,1)$-SDOP problem (cf.~[\refBOZ, \S6(2ii)]).

\smallskip
Case 2.2. $a_0=0$. Then $\Delta=y^2\Delta_0$, $\Delta_0=c_1-b_1^2$.
We proceed as in the beginning of \S\sectALGtoDOP.
Let
$$
     L^1=p(x)=p_0+p_1x, \qquad L^2=q_1 y+q(x)                    \eqno(\eqL)
$$
($p_0,\,p_1,\,q_1$ are constants).
Then the equation (\eqDkDj) takes the form $e_2y^2+e_1y+e_0=0$ where
$e_2=q_1\Delta_0'$.
The equation $e_2=0$ yields $c_1-b_1^2=\const$. Then we reduce
to the $(1,1)$-problem by the same variable change as in Case 2.1
 (cf.~[\refBOZ, \S6(2iii)]).

\smallskip
Case 3. $\deg_y\Gamma=0$. Then $\Omega$ contains a vertical strip $\Omega_0$.
We have $\deg_y a\le 1$.
If $\deg_y a=0$, then this is a solution of the $(1,\infty)$-SDOP problem
(see Figure~\figABCii). If $\deg_y a=1$, then $a$ vanishes at some point
$P\in\Omega_0$, thus $\g$ cannot be positive definite in $P$.
\qed\enddemo

%%%%%%%%%%%%%%%%%%%%%%%%%%%%%%%%%%%%%%%%%%%%%%%%%%%%%%%%%%%%%%%
%%%%%%%%%%%%%%%%%%%%%%%%%%%%%%%%%%%%%%%%%%%%%%%%%%%%%%%%%%%%%%%
%%%%%%%%%%%%%%%%%%%%%%%%%%%%%%%%%%%%%%%%%%%%%%%%%%%%%%%%%%%%%%%
%%%%%%%%%%%%%%%%%%%%%%%%%%%%%%%%%%%%%%%%%%%%%%%%%%%%%%%%%%%%%%%
%%%%%%%%%%%%%%%%%%%%%%%%%%%%%%%%%%%%%%%%%%%%%%%%%%%%%%%%%%%%%%%

\proclaim{ Theorem \thNonCInf }
Let $(\Omega,g,\rho)$ be a solution of the $(1,\infty)$-SDOP problem in $\bR^2$
with an unbounded $\Omega$.
Then one of the following cases occurs
up to a $(1,\infty)$-admissible change of variables.

\medskip\noindent
\roster
\item"(\modelUi)"
      $(\g,\Gamma)$ is as in Proposition~\propDegGaTwo(ii) with
      $\Omega=\{x>0\}\cap\{x^n>y^2\}$, $n\ge 1$, $c_{02}\le 0$,
      and $\rho$ is as follows:
\smallskip\noindent
      \hbox{\hskip-10pt\vbox{
      \roster
        \item"$\bullet$"
          if $n$ is odd, $\rho=(x^n-y^2)^{p-1} e^{-\lambda x}$, $\lambda>0$,
          $p>\max(0,\frac12-\frac1n)$;
        \smallskip
        \item"$\bullet$"
          if $n$ is even,
          $\rho=(x^{\frac{n}2}+y)^{p-1} (x^{\frac{n}2}-y)^{q-1} e^{-\lambda x}$,
          $\lambda>0$, $p+q>1-\frac2n$.
      \endroster
      }}
\smallskip\noindent
\item"(\modelSum)"
      $(\Omega,\g,\rho)$ is a product of two one-dimensional solutions, i.e.,
      $\Omega=\Omega_1\times\Omega_2$, $\g=\diag\big(\alpha g_1(x),\beta g_2(y)\big)$,
      $(\alpha,\beta)\ne(0,0)$, and
      $\rho=\rho_1(x)\rho_2(y)$ where each of $(\Omega_k,g_k,\rho_k)$ is
      one of the three one-dimensional solutions mentioned in Remark~\remDimOne\
      (which correspond to Hermite, Laguerre, and Jacobi polynomials);
      see also Remark~\remCorrig.
      
\endroster
\endproclaim

\demo{ Proof }
Let $w_{\min}$ be the minimal number $w$ such that $w\ge 1$ and
$(\Omega,g,\rho)$ is a solution of the $(1,w)$-SDOP problem.
If $w_{\min}=1$, the result can be derived from
the classification in [\refBOZ, \S\S5--6]. So, we assume that $w_{\min}>1$.
Let $\Delta=\deg\g$ and $\Gamma$ be the maximal boundary for $\g$
(see Definition~\defMaxBndry). Then $\partial\Omega\subset\{\Gamma=0\}$.
We have $\deg_y\Gamma\le\deg_y\Delta\le 2$ (see Figure~\figABCii).

\medskip
{\bf Case 1.} $\deg_y\Gamma=2$. Then $(\Gamma,\g)$ realizes one of Cases (i)--(v)
of Proposition~\propDegGaTwo. In Case (iii) with $n=0$ and in Cases (iv)--(v) 
we obtain (\modelSum) by Proposition~\propOPlus.
In Case (i), the same arguments as in the proof
of Theorem~\thCompInf\ show the absence of the exponential factor in $\rho$,
as well as in Case (iii) with $x_0=1$, $n>0$.

\smallskip
In Case (iii) of Proposition~\propDegGaTwo\ with $x_0=0$, $n>0$, the above arguments
(based on Proposition~\propMes) do not apply
literally, because $\Delta$ is no longer squarefree,
but they give the result being combined with [\refBOZ, Corollary~2.19].
Let us give however a self-contained proof of the unrealizability of this case.
We follow the algorithm from \S\sectALGtoDOP. Let $h=\log\rho$.
Rewrite $\g$ replacing $x$ by $-x$ and changing the sign:
$$
   \g = \left(\matrix x^2          & \tfrac12 nxy \lower5pt\hbox{\mathstrut}\\
                      \tfrac12 nxy & \tfrac14 n^2x^n - c_{02}\Gamma_2\endmatrix\right),
   \qquad \Gamma_2 = x^n - y^2,
   \qquad \Delta = \big(\tfrac14 n^2-c_{02}\big)x^2\Gamma_2.
$$
Let $L^i$ be as in (\eqL).
Then (\eqDkDj) reads $np_0 y- nxq(x) + 2x^2 q'(x) = 0$.
Hence $p_0=0$ and $q(x) = Cx^{n/2}$ ($C=0$ when $n$ is odd).
Plugging the obtained solution into (\eqDh) and integrating $h'_x$ and $h'_y$,
we obtain $\rho = C_1x^\alpha(x^n-y^2)^\beta$ when $n$ is odd, and
$\rho= C_1x^\alpha(x^{n/2}+y)^{\beta_1}(x^{n/2}-y)^{\beta_2}$ when $n$ is even,
with some constants $C_1$, $\alpha$, $\beta$, $\beta_1$, $\beta_2$.
Thus $\rho$ is not integrable on any component of
$\bR^2\setminus\{\Gamma=0\}$.

\smallskip
In Case (ii) of Proposition~\propDegGaTwo, we obtain (\modelUi). Indeed,
we have $(a,b,c)=(x,\frac12 ny,\frac14 n^2 x^{n-1}-c_{02}\Gamma)$ and
$\Delta=\Gamma_0\Gamma$ where $\Gamma_0=\frac14 n^2-c_{02}x$ and $\Gamma=x^n-y^2$.
It is clear that $n>0$ because if $n=0$, then $\Delta$ vanishes in each component
of $\{\Gamma=0\}$ which implies that $\g$ cannot be positive definite in $\Omega$.
The same argument leads to our choice of the real form of $\Gamma$ and to
the choice of $\Omega$ (up to the variable change $x\mapsto-x$).
The form of $\rho$ follows from Proposition~\propMes. Indeed, let $Q$ be as in
(\eqMes). Since $(\Omega,\g,\rho)$ is a solution of $(1,w)$-SDOP problem for some 
$w>2$, we deduce from (\eqMesIneq)
that $\deg_yQ=0$. But $(\Omega,\g,\rho)$ is also a solution for
$w=\frac12\deg c_0$ (thus $w=n$ if $c_{02}\ne 0$ and $w=n-1$ if $c_{02}=0$), and
we have $\deg_x\Delta=w+1$. Thus (\eqMesIneq) implies $\deg_xQ\le 1$.
The inequalities for
$p$, $q$, and $\lambda$ are the integrability conditions (see [\refBOZ, Remark~2.28]),
and $c_{02}\le 0$ is equivalent to $\Gamma_0|_\Omega>0$ and hence to
$\g|_\Omega>0$. Thus we obtain (\modelUi).

\medskip
{\bf Case 2.} $\deg_y\Gamma=1$. Then $(\Gamma,\g)$ realizes one of Cases (i)--(vi)
of Proposition~\propDegGaOne.

\smallskip
Case 2(i).
Then $k=1$ because otherwise $\Gamma$ is not $\g$-maximal.
We have $\Gamma=x\Gamma_1$, $\Gamma_1=x^ny-1$,
$\Delta=-x^2c_0(x)\Gamma_1$, thus $c_0\ne 0$.
As above, we follow \S\sectALGtoDOP\ to find $\rho$.
Let $L^i$ be as in (\eqL). Then (\eqDkDj) takes the form
$$
  \Big((np+q_1x)xc'_0 + (p_0+2np+2q_1x)nc_0\Big)y + nxqc_0+x^2(qc_0'-q'c_0)=0,
$$
Equating the coefficient of $y^0$ to zero, we obtain $q=C_1 x^nc_0$ 
(here and below $\alpha,\,\beta,\,C_1,\,C_2,\,C_3$ are some constants).
Equating the coefficient of $y^1$ to zero, we obtain two solutions: the first one
is $p_0=0$, $q_1=-np_1$ (then $c_0$ is an arbitrary function);
the second one is $c_0=C_2 (np+q_1x)/x^{2n+1}$. The first solution yields
$\rho=C_3 x^\alpha\Gamma_1^\beta$ which contradicts the integrability condition.
The second solution is irrelevant because it cannot be a nonzero polynomial.

\smallskip
Case 2(ii).
We have $\Gamma=xy-1$ and $\Delta=(xy+1-x^2c_0-b_{11}^2\Gamma+2b_{11})\Gamma$.
Then $c_0\ne 0$ (otherwise $w_{\min}=1$) and
$b_{11}\ne-1$ (otherwise $\g$ is as in 2(i), hence $\Gamma$ is not $\g$-maximal).
Hence
$\Gamma^2$ does not divide $\Delta$, thus we can apply
Proposition~\propMes\ and write $\rho$ in the form (\eqMes).
Then (\eqMesIneq) with $\ww=(1,1+\deg c_0)$
implies $\deg_x Q=0$ and (\eqMesIneq) with $\ww=(1,w)$, $w\gg0$, implies $\deg_y Q=0$.
Hence $\rho$ is not integrable.

\smallskip
Case 2(iii).
We have $(a,b,c;\Gamma)=(a_0,b_1y,c_2y^2+c_1y;y)$, $\Delta=(a_0c_2-b_1^2)y^2+a_0c_1y$.
If $\deg c_1<2$, then $w_{\min}=1$, thus $\deg c_1\ge 2$.
Suppose that $\deg a\ge 1$.
Then $\deg_x ac\ge 3>\deg_x b^2$ whence $\deg_x\Delta=\deg_x ac$
and we obtain a contradiction with Proposition~\propMesII\ for $w=\deg c_1$.
Thus $0\ne a=\const$ and we may set $a=1$.

We have $ac_1\ne 0$, hence $y$ does not divide $\Delta$.
Then Proposition~\propMes\ with $\ww=(1,w)$,
$w\gg0$, implies that $\rho=y^p e^h$ with $\deg_y h=2-\deg_y\Delta$.
If $\deg_y\Delta=2$, this contradicts the integrability condition,
hence $\deg_y\Delta=1$. Then $c_2=b_1^2$, in particular $b_1=\const$,
thus $(a,b,c)=(1, \beta y, \beta^2 y^2 + c_1 y)$, $\beta\in\bR$.
Let $h=h_1(x)y+h_0(x)$.
Then $\deg_y L^1=0$ for $L^1=ah'_x+bh'_y=h'_1y+h'_0+\beta yh_1$, hence
$h'_1=\beta h_1$. If $\beta\ne 0$, we obtain $h_1=0$ (since $h_1$ is
a polynomial). If $\beta=0$, the condition $\deg_{(1,w)}L^2\le w$ for
$L^2 = b h'_x + c h'_y = yc_1 h_1$ again implies $h_1=0$ because
$\deg c_1>0$. Thus $\deg_y h=0$
which contradicts the integrability condition.

\smallskip
Case 2(iv).
Then $(a,b,c;\Gamma)=(x\tilde a_0,\,b_{11}xy,\,c_{2}y^2+c_1y;\,xy)$,
$\Delta=(\tilde a_0c_{2}-b_{11}x)xy^2+ac_1y$. As in Case~2(iii),
we obtain $\deg c_1\ge 2$,
$\tilde a_0=1$, $y$ does not divide $\Delta$, hence (see Remark~\remMesBis)
$\deg_y\Delta=1$, that is $c_{2}=b_{11}^2x$.
Thus we arrive to $(a,b,c)=(x,\beta xy,\beta y^2+c_1y)$.
Then $\rho=x^py^qe^h$ with $h=h_1y+h_0$ where $h_k$ are rational functions of $x$.
The rest of the proof is as in Case~2(iii).

\smallskip
Case 2(v).
We obtain (\modelSum) by Proposition~\propOPlus.

Case 2(vi).
$a=0$ is impossible for a positive definite $\g$.
 
\medskip
{\bf Case 3.} $\deg_y\Gamma=0$, i.e., $\Gamma=\Gamma(x)$ is a polynomial in $x$ only.
Then $\Omega=I\times\bR$ for a finite or infinite (from either side) interval $I$.
Further, $\rho$ is of the form (\eqMes) (see Remark~\remMes)
where $\deg_y Q\ge 2$ by the integrability condition,
hence $\deg_y Q=2$ and $\deg_y\Delta=0$ by (\eqMesIneq);
we write $Q=h_2y^2+h_1y+h_0$ where $h_k$ are rational functions of $x$ and $h_2\ne0$.
Then (see (\eqDh))
$$
  \matrix
  \deg_\ww L^1=1\quad\text{for}\quad L^1 = a(h_2'y^2 + h_1'y + h_0')+b(2h_2y + h_1),\\
  \deg_\ww L^2=w\quad\text{for}\quad L^2 = b(h_2'y^2 + h_1'y + h_0')+c(2h_2y + h_1).    
  \endmatrix
                                                                        \eqno(\eqDhh)
$$
Since $\deg_y\Gamma=0$, (\eqParam) implies that $\Gamma$ divides $a$ and $b$,
thus $\deg_x\Gamma\le 2$.
Since $\Omega\subset\{\Gamma=0\}$, Proposition~\propMesII\
implies
$$
    \deg_x\Delta  \le \deg_x\Gamma - 1
     + \max\big(\,\lfloor w_{\min}\rfloor + \deg_x b,\,1+\deg_x c\,\big).   \eqno(\eqM)
$$

\smallskip
Case 3.1. 
$\deg_x\Gamma=2$. Then, up to rescaling,
$a=\Gamma$ and $b=b_0=\tilde b(x)\Gamma$ for some polynomial $\tilde b$. Then
the variable change $y\mapsto y-p(x)$ with $p'=\tilde b$ (see Example~\exaChangeVarII)
makes $b=0$. Then $\deg_y\Delta=0$ implies $c=c_0$.
Then (\eqDhh) yields $h_2'=0$ and $2h_2c_0=\const$.
Since $h_2\ne0$, we conclude that $c_0=\const$, hence we obtain (\modelSum).

\smallskip
Case 3.2. $\Gamma=x$.
Then $a=x\tilde a$ and 
$b=x\tilde b=x(\beta y+\tilde b_0(x))$ ($\beta=b_{11}$) by (\eqParam),
hence the coefficient of $y^2$ in $\Delta$ is $ac_{2}-x^2\beta$.
Since $\deg_y\Delta=0$, $a\ne0$, and $c_{2}=\const$, we have either $a=a_{20}x^2$
(then we may assume $a_{20}=1$) or $\beta=0$.

\smallskip
Case 3.2.1. $a=x^2$.
We have $\Delta=x^2 d(x)$ for
some polynomial $d=d(x)$, so we may write $c=\tilde b^2+d$.
If $\deg\tilde b_0\le 1$, then either (\eqM) fails or $w_{\min}=1$.
Thus we assume that $\deg\tilde b_0\ge 2$.
The change of variables $y\mapsto y+\lambda x^n$
transforms $b$ into $b+\lambda(n-\beta)x^{n+1}$ (see Example~\exaChangeVarII).
Hence we may kill all coefficients of $b_0$ unless $\beta=n\in\bN$
in which case we kill all of them except $b_{n+1,0}$.
Since $\deg\tilde b_0\ge2$,
we then assume that $\tilde b=ny+\alpha x^n$, $n\ge 2$.

Then (see (\eqDhh))
the coefficient of $y^2$ in $L^1$ is $x^2h_2'+2nxh_2$, hence $h_2'=-2nxh_2$ whence
$h_2=Cx^{-2n}$ which contradict the facts that $h_2\ne0$ and 
$h_2$ is a rational function with denominator at most $x$ (see [\refBOZ, Prop.~2.15]).

\smallskip
Case 3.2.2. $\beta=0$. Then the condition $\deg_y\Delta=0$ implies $c_2=c_1=0$,
i.e., $\g$ does not depend on $y$.
By the change of variables $y\mapsto y+p(x)$ we may achieve that
$\deg_x\tilde b<\deg_x\tilde a$.
If $\deg\tilde a=0$, then $b=0$ and the proof is
the same as in Case~3.1. If $\deg\tilde a=1$ and $b\ne0$,
we obtain a contradiction with (\eqM).

\smallskip
Case 3.3. $\Gamma$ is constant. Then $\Omega=\bR^2$. Since $a|_\Omega>0$,
we have up to translation $a=x^2+1$ or $a=1$. 

\smallskip
Case 3.3.1. $a=x^2+1$. Recall that $\deg_y\Delta=0$,
hence $(x^2+1)c_2=b_1^2$ and $ac_1=2b_1b_0$ whence $b_1=c_2=c_1=0$.
By a variable change $y\mapsto y-p(x)$ we may achieve that $\deg b\le 1$
(see Example~\exaChangeVarII) and we obtain a contradiction with (\eqM).

\smallskip
Case 3.3.2. $a=1$.
If $b_1=0$, the proof is as in Case~3.1. If $b_1\ne 0$,
then (\eqDhh) for $L^1$ implies $h_2'=-b_1 h_2$ which is impossible
for nonzero polynomials.
\qed\enddemo

\medskip\noindent
{\bf Remark \remCorrig.} There is a small mistake in [\refBOZ, \S4.2].
It is erroneously claimed that each cometric solution on the square $[-1,1]^2$
is proportional to $\diag(1-x^2,1-y^2)$. The correct answer is
$\g=\diag\big(\alpha(1-x^2),\beta(1-y^2)\big)$ with arbitrary positive $\alpha$
and $\beta$. All the corresponding riemannian metrics $(g_{ij})=\g^{-1}$ are
pairwise non-isometric.

%%%%%%%%%%%%%%%%%%%%%%%%%%%%%%%%%%%%%%%%%%%%%%%%%%%%%%%%%%%%%%%%%
%%%%%%%%%%%%%%%%%%%%%%%%%%%%%%%%%%%%%%%%%%%%%%%%%%%%%%%%%%%%%%%%%
%%%%%%%%%%%%%%%%%%%%%%%%%%%%%%%%%%%%%%%%%%%%%%%%%%%%%%%%%%%%%%%%%
%%%%%%%%%%%%%%%%%%%%%%%%%%%%%%%%%%%%%%%%%%%%%%%%%%%%%%%%%%%%%%%%%
%%%%%%%%%%%%%%%%%%%%%%%%%%%%%%%%%%%%%%%%%%%%%%%%%%%%%%%%%%%%%%%%%
%%%%%%%%%%%%%%%%%%%%%%%%%%%%%%%%%%%%%%%%%%%%%%%%%%%%%%%%%%%%%%%%%

\head\sectRealiz. Realization of solutions as images of the Laplace operator
\endhead

Let $M_1$ and $M_2$ be smooth manifolds and $\Phi:M_1\to M_2$ be
a smooth mapping which is a submersion at a generic point of $M_1$.
Let $\cL_1$ and $\cL_2$ be differential operators on $M_1$ and on
$\Omega=\Phi(M_1)$ respectively. We say that $\cL_2$ is
{\it the image of $\cL_1$ through} $\Phi$ if
$\cL_2(f) = \cL_1(f\circ\Phi)$.
Notice that the image of $\cL_1$ through $\Phi$ may or may not exist,
moreover, it does not exist for generic $\cL$ and $\Phi$ unless $\Phi$ is
injective. However it does exist when $\cL_1$ and $\Phi$ are invariant under an action
of a group $G$ on $M_1$, and then $\Phi$ identifies $\Omega$ with the
orbit space $M_1/G$.

For example, for half-integer $p$ and $q$,
the Jacobi operator $J_{p,q}$ on the interval $(-1,1)$ (see Remark~\remDimOne) is the
image of the Laplace operator $\De_{\bS^n}$
on the sphere $\bS^n$, $n=2p+2q-1$,
through the mapping
$$
  (x_1,\dots,x_{2p},y_1,\dots,y_{2q})\mapsto x_1^2+\dots+x_{2p}^2
$$
(see [\refBOZ, \S2.1], the end of [\refBOZ, \S4.1], and references therein).

In [\refBOZ, \S4]
similar interpretations are given for many
values of parameters of each compact solution of $(1,1)$-DOP problem.
They are interpreted as images of the Laplace (Casimir)
operators on $\bS^n$, $\bR^n$, $SO(n)$, or $SU(n)$.

Each quotient of $\bS^2$ or $\bR^2$ by a reflection group
can be identified with a compact domain $\Omega$,
$\partial\Omega\subset\{\Gamma(x,y)=0\}\subset\bR^2$, so that
$(\Omega,\LL,\Gamma^{-1/2}dx)$ is
a solution of $(1,w)$-DOP problem where $\LL$ is the image
of $\De_{\bS^2}$ or $\De_{\bR^2}$ through the quotient map.
Under this identification, the vertices of the fundamental polygon of the group
action are in bijection with singular points of $\partial\Omega$ so that
the angle $\pi/n$ corresponds to the singularity $A_{n-1}$
(given by $u^2=v^n$ in some local coordinates).
Explicit formulas for the mappings $\bS^2\to\bR^2$ and $\bR^2\to\bR^2$ which
realize the operators $\LL$ as images of the Laplace operator are given in
[\refBOZ] and [\refBB]; see detailed references in Table~\tabQuo. 

\midinsert
\vbox{
\def\TL#1#2#3#4#5{%
 \noindent\hskip0pt
 \hbox to 6mm{\hskip3pt#1\hfill}
 \hbox to 12mm{#2\hfill}
 \hbox to 76mm{#3\hfill}
 \hbox to 3.5mm{#4\hfill}
 \hbox to 20mm{#5\hfill}}
\hrule
\vskip5pt
\TL{}{Angles}{Boundary of $\Omega$}{$w$}{Reference}
\vskip3pt
\hrule
\vskip3pt
\TL{}{2,2,2,2}{Rectangle (see also Remark~\remCorrig)}{1}{[\refBOZ,\S4.2]}

\TL{\lower6pt\hbox{$\bR^2$}}
     {2,4,4}{Parabola with two tangents\hfill}{1}{[\refBOZ,\S4.7]}
\vskip-3pt

\TL{}{3,3,3}{Deltoid: $x=2\cos\theta+\cos2\theta$, $y=2\sin\theta-\sin2\theta$}{1}
                    {[\refBOZ,\S4.12]}

\TL{}{2,3,6}{Cubic $y^2=x^3$ with a cubically tangent parabola}{2}
                    {[\refBOZ, end of \S4.12]}
\vskip3pt
\hrule
\vskip3pt
\TL{}{---}{Circle}                {1}{[\refBOZ,\S4.3]}

\TL{}{2,2}{Coaxial parabolas: $y^2=(1-x^2)^2$}{1}{[\refBOZ,\S4.5]}

\TL{}{$n,n$}{$y^2=(1-x^2)^n$, $n\ge 2$}{$n$}{[\refBB,\S6]: $\Omega_1^{(n)}$}

\TL{}{2,2,2}{Triangle}{1}{[\refBOZ,\S4.4]}

\TL{\lower6pt\hbox{$\bS^2$}}
     {2,2,3}{$(y^2-x^3)(x-1)=0$}{1}{[\refBOZ,\S4.9]}
\vskip-3pt

\TL{}{2,2,4}{$y(y-x^2)(x-1)=0$}{1}{[\refBOZ,\S4.6]}

\TL{}{2,2,$n$}{$(y^2-x^n)(x-1)=0$, $n\ge2$}{$n$}
                                      {[\refBB,\S6]: $\Omega_3^{(n)}\!\!,\Omega_6^{(n)}$}

\TL{}{2,3,3}{Swallow tail: $\text{discrim}_t(t^4-t^2+xt+y)=0$}{1}{[\refBOZ,\S4.11]}

\TL{}{2,3,4}{Cubic $y^2=x^3$ with a tangent line}{1}{[\refBOZ,\S4.10]}

\TL{}{2,3,5} {Dodecahedral quotient, see Theorem~\thCompTwo(\modelBi)}
                    {2}{[\refBB,\S8]: $\Omega_{21}$}
\vskip3pt
\hrule
}   % end of \vbox
\botcaption{Table \tabQuo}
  Quotients of $\bR^2$ and $\bS^2$ by reflection groups; $a,b,\dots$ in the ``Angles''
  cell means: the angles of the fundamental domain are 
  $\frac{\pi}{a},\frac{\pi}{b},\dots$
\endcaption
\endinsert

If $\LL$ lifts to $\De_{\bS^2}$ or $\De_{\bR^2}$,
then the corresponding metric is of constant non-negative curvature.
Moreover, the classification obtained in [\refBOZ] implies that the
curvature is constant and non-negative for all
solutions with $\det(\g)$ of maximal degree.
For maximal degree solutions, Lev Soukhanov [\refSouArx], [\refSouAFST]
proved the constancy of curvature (but not its non-negativity)
not using the classification. He also proved a
generalization of this result to an arbitrary dimension.

The curvature in solution (\modelBiv) in Theorem~\thCompInf\
is constant if and only if $m=n$ and $c_{02}=-n^2$. In (\modelBv) it
is constant if and only if $c_{02}=-\frac14 n^2$.
If $m\ne n$ in (\modelBiv), the curvature cannot be constant
because there are no biangles of constant curvature with different angles.
For $m=n$ in (\modelBiv) and for (\modelBv), the curvature is
$$
     \frac{-\lambda n^2\big(2(n^2+\alpha)x^k+\alpha\big)}
            {\big((n^2+\alpha)x^k-\alpha\big)^2}, \qquad
     (k,\lambda,\alpha) = \cases
                             (2,2,c_{02}),        &\text{(\modelBiv), $m=n$},\\
                             (1,\frac12,4c_{02}), &\text{(\modelBv)}, \endcases
$$
which is evidently non-constant when $\alpha\ne-n^2$ (cf. [\refBOZ, \S4.5]).
The similarity of the formulas for the curvature is not occasional:
up to adjusting the constants, (\modelBv) is the image of (\modelBiv) through
$(x,y)\mapsto(x^2,y)$.

According to our classification, the only bounded domains which admit solutions
but are not covered by
Table~\tabQuo\ are those in Theorem~\thCompInf(\modelBiv) for $m\ne n$.
The simplest one is the nodal cubic $y^2=x^2-x^3$ which corresponds to $(m,n)=(1,2)$.
This solution is realized in [\refBOZ, \S4.8] as the image of $\De_{\bS^{3c}}$,
$c=1,2,4,8$.
The realization for $c=1$ immediately extends to any $(m,n)$ as follows
(it could be interesting to do the same for $c=2,4,8$).
We consider $\bS^3$ as the unit sphere in $\bC^2$ with coordinates $(z_1,z_2)$.
Then the image of $\frac14\De_{\bS^3}$ through
$$
   (z_1,z_2)\mapsto(X,Y)=\big(\,|z_1|^2,\, \Re(z_1^n\bar z_2^m)\,\big)
$$
is an operator on the domain bounded by the curve
$(1-X)^m X^n - Y^2 = 0$. Its image through the affine change
$(X,Y)\mapsto(x,y)=(\,2X-1,\,2^{(m+n)/2}Y\,)$
is the operator corresponding to (\modelBiv) with $p=q=\frac12$ and
$c_{02}=-\frac14(m+n)^2$.

\enddocument